\documentclass[10pt, a4paper]{article}
\usepackage{graphics}
\usepackage{latexsym}
\usepackage{amssymb}
\usepackage{amsmath}
\usepackage{color}
\usepackage{comment}

\newcommand{\beq}{\begin{equation}}
\newcommand{\eeq}{\end{equation}}
\newcommand{\beqas}{\begin{eqnarray*}}
\newcommand{\eeqas}{\end{eqnarray*}}
\newcommand{\ep}{\varepsilon}
\newcommand{\eps}{\epsilon}
\newcommand{\ue}{u^{\varepsilon}}

\newcommand{\wc}{\rightharpoonup}

\newcommand{\om}{{\omega}}

\newcommand{\bY}{{\bf Y}}
\newcommand{\by}{\boldsymbol{y}}
\newcommand{\bV}{{\bf V}}
\newcommand{\bX}{{\bf X}}
\newcommand{\bcX}{{\boldsymbol{\cal X}}}
\newcommand{\bv}{\boldsymbol{v}}
\newcommand{\bu}{\boldsymbol{u}}

\newcommand{\bH}{{\bf H}}
\newcommand{\cH}{{\cal H}}
\newcommand{\bcH}{{\boldsymbol{\cal H}}}

\newcommand{\bz}{{\boldsymbol{z}}}

\newcommand{\bep}{\boldsymbol{\epsilon}}
\newcommand{\ubV}{{\underline{\bV}}}
\newcommand{\ubH}{{\underline{\bH}}}
\newcommand{\ubcH}{{\underline{\bcH}}}
\newcommand{\ubX}{{\underline{\bX}}}
\newcommand{\ubcX}{{\underline{\bcX}}}

\newcommand{\ubVL}{{\underline{\bV}_\Lambda}}
\newcommand{\ubHL}{\underline{\bH}_\Lambda}
\newcommand{\ubcHL}{{\underline{\bcH}_\Lambda}}

\newcommand{\ubcXL}{\underline{\bcX}_\Lambda}
\newcommand{\buL}{{\bu}_\Lambda}
\newcommand{\fb}{\frak b}
\newcommand{\fa}{\frak a}
\newcommand{\fd}{\frak d}

\newcommand{\be}{\begin{equation}}
\newcommand{\ee}{\end{equation}}
\newcommand{\bea}{\begin{eqnarray}}
\newcommand{\eea}{\end{eqnarray}}
\newcommand{\beas}{\begin{eqnarray*}}
\newcommand{\eeas}{\end{eqnarray*}}
\newcommand{\ds}{\displaystyle}

\def\IN{\mathbb{N}}

\def\IP{\mathbb{P}}
\def\IR{\mathbb{R}}

\def\cF{{\cal F}}

\def\bproof{{\it Proof\ \ }}
\def\eproof{{\hfill$\Box$}}

\newtheorem{theorem}{Theorem}[section]
\newtheorem{lemma}[theorem]{Lemma}
\newtheorem{assumption}[theorem]{Assumption}
\newtheorem{proposition}[theorem]{Proposition}
\newtheorem{definition}[theorem]{Definition}
\newtheorem{remark}[theorem]{Remark}

\definecolor{darkgreen}{rgb}{0,0.7,0}

\numberwithin{equation}{section}
\title
{
Polynomial approximations of a class 
of  stochastic multiscale elasticity problems
\thanks{%
This research was supported by the AcRF Tier 1 grant RG69/10, the Singapore A*Star SERC grant 122-PSF-0007, the AcRF Tier 2 grant MOE 2013-T2-1-095, ARC 44/13 and a graduate scholarship from Nanyang Technological University. 
}
}
\author{
Viet Ha Hoang
\thanks{
Division of Mathematical Sciences, 
School of Physical and Mathematical Sciences, 
Nanyang Technological University, Singapore, 637371
}
, 
Thanh Chung Nguyen
\thanks{
Division of Mathematical Sciences, 
School of Physical and Mathematical Sciences, 
Nanyang Technological University, Singapore, 637371; Permanent address: Department of Mathematics, Quang Binh Universitiy, 312 Ly Thuong Kiet, Dong Hoi, Vietnam.
}
and
Bingxing Xia
\thanks{
Division of Mathematical Sciences, 
School of Physical and Mathematical Sciences, 
Nanyang Technological University, Singapore, 637371;
}
}
\date{}
\begin{document}
\setcounter{page}{1}
\maketitle
\begin{abstract}
We consider a class of elasticity equations in $\IR^d$ whose elastic moduli depend on $n$ separated microscopic  scales. The moduli are random and expressed as a linear expansion of a countable sequence of  random variables which are independently and identically uniformly distributed in a compact interval. The multiscale displacement problem, the multiscale Hellinger-Reissner mixed problem that allows for computing the stress directly, and the multiscale mixed problem with a penalty term for nearly incompressible isotropic materials are considered. The stochastic problems are studied via deterministic problems that  depend on a countable number of real parameters which represent the probabilistic law of the  stochastic equations. We study the multiscale homogenized problems that contain all the macroscopic and microscopic information. The solutions of these multiscale homogenized problems are written as generalized polynomial chaos (gpc) expansions. We approximate these solutions by semidiscrete Galerkin approximating problems that project into the spaces of functions with only a finite number of $N$ gpc modes. Assuming summability properties for the coefficients of the elastic moduli's expansion, we deduce  bounds and summability properties for the solutions' gpc expansion coefficients. These bounds imply explicit rates of convergence in terms of $N$  when the gpc modes used for the Galerkin approximation are chosen to correspond to the best $N$ terms in the gpc expansion. For the mixed problem with a penalty term for nearly incompressible materials, we show that the rate of convergence for the best $N$ term approximation is independent of the Lam\'e constants' ratio when it goes to $\infty$.  Correctors for the homogenization problem are deduced. From these we establish correctors for the solutions of the parametric multiscale problems in terms of the  semidiscrete Galerkin approximations. For two scale problems, an explicit homogenization rate which is uniform with respect to the parameters is deduced. Together with the best $N$ term rate, it provides an explicit convergence rate for the correctors of the parametric multiscale problems. For nearly incompressible materials, we obtain a homogenization rate that is independent of the ratio of the Lam\'e constants, so that the error for the corrector  is also independent of this ratio.  
\end{abstract}

\section{Introduction}
We consider a multiscale elasticity problem in $\IR^d$ whose elastic tensor is random and is a linear combination of a sequence of  random variables which are independently and uniformly distributed in a compact interval. The elastic tensor depends on $n$ separable microscopic scales and is periodic with respect to each of these scales. We consider  the multiscale displacement problem, the multiscale Hellinger-Reissner mixed problem that allows for computing the stress tensor directly, and the multiscale mixed problem with a penalty term for nearly incompressible isotropic materials.  We study the random problems via  deterministic ones whose elastic tensor depends on an infinite sequence of real parameters. The space of parameter sequences is equipped with a probability measure that is the law of the sequence of random variables that the random elastic tensor depends on. Thus the deterministic parametric multiscale solution is the law of the solution of the stochastic multiscale problem. Solving this parametric equation, we obtain  statistic properties of the stochastic multiscale solution. 

For multiscale problems, a direct numerical procedure is prohibitively expensive. The problem is approximated by the homogenization limit when all the microscopic scales converge to 0 (\cite{BLP},\cite{JKO}). To obtain the microscopic information, a part from the solution of the homogenized problem, we need also the scale interacting terms. We therefore apply the multiscale convergence (see \cite{Nguetseng}, \cite{Allaire} and \cite{AllaireBriane}) to obtain the  multiscale homogenized problem that contains all the macroscopic and microscopic information. The problem is posed in a high dimensional tensorized domain: if the original multiscale problem is posed in $\IR^d$ and depends on $n$ microscopic scales and one macroscopic scale, this problem is posed in $\IR^{(n+1)d}$. However, as demonstrated in Hoang and Schwab \cite{HSelliptic} for multiscale elliptic equations and in Xia and Hoang \cite{XHelasticity} for multiscale elasticity equations, the sparse tensor finite element approach is capable of solving these high dimensional multiscale homogenized problems with an essentially optimal complexity which is essentially equal to that for solving a problem in $\IR^d$. Though we do not address finite element approximation  in the paper, this is the motivation for us to consider polynomial approximations for the solutions of  the stochastic/parametric multiscale homogenized problems.  We
write their solutions  in terms of a generalized polynomial chaos (gpc) expansion with respect to a system of multivariate polynomials which forms an orthonormal basis for the $L^2$ Lebesgue space of the parameter sequences. Following  Cohen et al. \cite{CDS1} and \cite{CDS2} and other related papers (\cite{HSmultirandom}, \cite{HSwave}, \cite{HSparabolic}), we study the best $N$ term approximation for the gpc expansion of the high dimensional multiscale homogenized problems. When assuming summability for the coefficients of the elastic moduli's expansion, we deduce summability properties for the coefficients of the gpc expansion. From this, an explicit error estimate for the semidiscrete Galerkin approximating problem which projects into the spaces of functions with only $N$ gpc modes can be deduced when these modes are chosen to correspond to the best $N$ terms in the gpc expansion. In many cases, this rate is superior over the Monte Carlo $N^{-1/2}$ rate.  

To approximate the solutions of the multiscale problems, we derive  correctors from the solution of the semidiscrete Galerkin problems. For two scales, an explicit homogenization error is available. To employ this error for the parametric problem, we prove that it is uniform with respect to the parameter sequences. From this a corrector for the multiscale problem in the mean square norm with respect to the parameter space is deduced. The rate of convergence is the sum of the best $N$ term semidiscrete Galerkin rate and the uniform homogenization rate. For more than two scale problems, an explicit homogenization rate is not available. However, we can construct a corrector for the parametric homogenization problem. From this, a corrector for the parametric multiscale problem in the mean square norm using the solution of the semidiscrete Galerkin problem is deduced, without an explicit error.  

For multiscale elliptic problems, the framework has been applied in \cite{HSmultirandom} where the complex method to deduce bounds for the gpc coefficients is employed. Here we employ the real method developed in \cite{CDS1}. However, the main contributions of this paper are the studies of the stochastic/parametric high dimensional multiscale homogenized problems for the multiscale mixed Hellinger-Reissner formula and for the multiscale mixed problem with a penalty term for nearly incompressible materials. We adapt the approach in \cite{XHelasticityrandom} for single scale macroscopic equations to deduce the best $N$ term approximations.  For the mixed problems, the inverse of the elastic tensor depends on the random variables in the expansion nonlinearly. To employ the linear dependence, we formulate the equivalent mixed problems that use the elastic tensor instead of its inverse. For the Hellinger-Reissner problem, showing the well-posedness of the semidiscrete Galerkin problem of the equivalent form, and deducing bounds for the coefficients of the gpc expansion require careful manipulation of various inf-sup conditions. For the mixed problem with a penalty term, we obtain a best $N$ term approximation whose rate is independent of the ratio of the Lam\'e constants. We also obtain a homogenization rate of convergence for nearly incompressible materials that is independent of this ratio. To the best of our knowledge, this homogenization rate is new. Together with the best $N$ term approximation, we construct a corrector for the stochastic/parametric multiscale nearly incompressible problem with an error independent of the ratio of the Lam\'e constants. 

The paper is organized as follows. In the next section, we set up the multiscale problem, define the multiscale and the random structures of the elastic tensors. We formulate  the multiscale random displacement problem,  the  multiscale mixed Hellinger-Reissner problems and the multiscale mixed problem with a penalty term for nearly incompressible materials. Section 3 studies the deterministic parametric problems. We show that their solutions are measurable with respect to  the $\sigma$-algebra of the parameter space. Therefore the solution of the original stochastic problem can be deduced from the parametric solution by inserting the random sequences into the place of the parameters. We establish the multiscale homogenized problems in this section. The next three sections are devoted to approximating the deterministic parametric high dimensional multiscale homogenized problems. Approximations for the displacement problem is considered in Section 4. We first write the high dimensional solution as a gpc expansion. We then consider the semidiscrete Galerkin problem that projects into a subspace of functions with only $N$ fixed gpc modes. To get an explicit rate of convergence, we deduce bounds for the coefficient functions of the gpc expansion. From these bounds, we derive the rate of convergence  when the gpc modes are chosen to  correspond to the best $N$ terms in the gpc expansion. We consider approximation of the mixed high dimensional multiscale homogenized problem for the Hellinger-Reissner setting and its equivalent form using the elastic tensor (but not its inverse) in Section 5.  Employing standard estimates for saddle point problems, we deduce bounds for the coefficient functions of the gpc expansions of the displacement vector and the stress tensor. We then deduce an explicit convergence rate for the semidiscrete Galerkin problem when the finite number of gpc modes is chosen corresponding to the best $N$ terms in the gpc expansion. Approximation for the nearly incompressible problem is studied in Section 6. We consider both the mixed problem with a penalty term and the equivalent one using the Lam\'e constant $\lambda$ instead of $1/\lambda$. We deduce bounds for the coefficients of the gpc expansion, from which we show that the best $N$ term approximation achieves a rate of convergence that is independent of the ratio of the Lam\'e constants when $\lambda$ goes to $\infty$. In Section 7, we use the semidiscrete Galerkin approximations to deduce correctors for the parametric multiscale problems. For two-scale problems, we derive a homogenization rate  which is uniform with respect to the parameters. From this a corrector for the parametric multiscale problem is deduced in the mean square norm with respect to the parameter space with an explicit error that is the sum of the uniform homogenization error and the best $N$ term approximation rate. For nearly incompressible materials, we prove a homogenization rate that is  independent of the ratio of the Lam\'e constants. The error for the correctors of the stochastic/parametric two-scale problem is thus independent of this ratio.  For problems that depend on more than one microscopic scales, we derive a corrector for the parametric homogenization problem, without a rate of convergence, which implies a corrector in the mean square norm with respect to the parameter space for the solution of the parametric multiscale problem, without an explicit rate of convergence. 

Throughout this paper, repeated indices indicate summation. Notation $\nabla$ without an explicit variable denotes the gradient with respect to $x$. Similarly, $\eps$ without an explicit variable denotes the elastic strain tensor with respect to $x$.  We denote by $:$ the inner product in $\IR^{d\times d}_{sym}$. For a sequence of integers $\nu=(\nu_1,\nu_2,\ldots)$ with only a finite number of terms being non zero, we denote by $\nu!=\nu_1!\nu_2!\ldots$; and for a sequence of real numbers $d=(d_1,d_2,\ldots)$, we denote by $d^\nu=d_1^{\nu_1}d_2^{\nu_2}\ldots$. The notation $\#$ indicates spaces of periodic functions, in particular $H^k_\#(Y)$ denotes Sobolev spaces of periodic functions, and $C^k_\#(Y)$ denotes the space of $k$ time differentiable periodic functions. 

\section{Setting-up of the problems}
Let $D\subset\IR^d$ be a bounded domain. Let $Y$ be the unit cube $(0,1)^d\subset\IR^d$ and $Y_1,\ldots,Y_n$ be $n$ copies of $Y$. We denote by $\bY=Y_1\times\ldots\times Y_n$ and $\by=(y_1,\ldots,y_n)\in\bY$. Let $(\Omega,\Sigma,\IP)$ be a probability space where $\Sigma$ is the sigma algebra and $\IP$ is the probability measure. The elastic tensor is a random function $a(\omega,x,y_1,\ldots,y_n):\Omega\times D\times \bY\to L^\infty(\Omega,C(D,C_\#(Y_1\times\ldots\times Y_n)))^{d^4}$
where $C_\#(Y_1\times\ldots\times Y_n)=C_\#(\bY)$ denotes the space of continuous functions that are periodic with respect to $y_i$ with the period $Y_i$ for $i=1,\ldots,n$.  The tensor function $a$ is assumed to be symmetric: for $i,j,k,l=1,\ldots,d$
\[
a_{ijkl} =a_{ijlk}=a_{klij}.
\]
The random tensor $a$ satisfies the coerciveness and boundedness condition
\be
a(\omega;x,\by)\xi:\xi\ge \alpha|\xi|^2,\ \ a(\omega;x,\by)\xi:\zeta\le \beta|\xi||\zeta|\ \ \forall\,\xi,\zeta\in\IR^{d\times d}_{sym}
\label{eq:boundednesscoerciveness}
\ee
where the constants $\alpha>0,\beta>0$ are independent of $\omega\in \Omega$, $x\in D$ and $\by\in \bY$. 
To define the $n$ microscopic scales  that the multiscale elasticity problem depends on, let $\ep>0$ be a small quantity and $\ep_i$ ($i=1,\ldots,n$) be $n$ positive functions of $\ep$ which satisfy the scale separation assumption
\[
\lim_{\ep\to 0}{\ep_{i+1}(\ep)\over\ep_i(\ep)}=0,\ \ i=1,\ldots,n-1.
\]
Without loss of generality, we assume that $\ep_1(\ep)=\ep$. We define the random multiscale elastic moduli as
\[
a^{\ep}(\omega;x)=a\left(\omega; x,{x\over\ep_1},\ldots,{x\over\ep_n}\right).
\]
To define the probability distribution of the random elastic moduli $a$, we assume that there are independent random variables $z_m:\Omega\to [-1,1]$ which are uniformly distributed in $[-1,1]$, and there are fourth order symmetric tensor functions $\psi_m:D\times\bY\to \IR^{d^4}$ so that
\be
a(\omega;x,\by)=\bar a(x,\by)+\sum_{m=1}^\infty z_m(\omega)\psi_m(x,\by).
\label{eq:randommoduli}
\ee
The fourth order tensor function $\bar a(x,\by)$ is the mean value of $a$. It is symmetric  and satisfies 
\be
\bar a(x,\by)\xi:\xi\ge\alpha_0|\xi|^2,\ \ |\bar a(x,\by)\xi:\zeta|\le\beta_0|\xi||\zeta|
\label{eq:abar}
\ee
for $\alpha_0>0,\beta_0>0$ and all $\xi,\zeta\in \IR^{d\times d}_{sym}$.  
For the uniform coerciveness and boundedness \eqref{eq:boundednesscoerciveness}, we assume further that there are positive constants $\beta_m$ such that for all $\xi,\zeta\in \IR^{d\times d}_{sym}$
\be
|\psi_m(x,\by)\xi:\zeta |\le \beta_m|\xi||\zeta|
\label{eq:beta}
\ee
and that 
\be
\sum_{m=1}^\infty\beta_m\le {\kappa\over1+\kappa}\alpha_0
\label{eq:alpha}
\ee
for $\kappa>0$. We can then take the constants $\alpha$ and $\beta$ in \eqref{eq:boundednesscoerciveness} as 
\be
\alpha={\alpha_0\over 1+\kappa},\ \ \beta=\beta_0+{\kappa\over1+\kappa}\alpha_0.
\label{eq:ab1}
\ee
Let $\Gamma$ be a subset of the boundary $\partial D$. Let $H^1_\Gamma(D)$ be the subspace of $H^1(D)$ of functions with zero trace on $\Gamma$.  
Let $V=H^1_\Gamma(D)^d$. We denote by $(\cdot,\cdot)$ the inner product in $L^2(D)^d$, extended to the dual pairing relation $\langle\cdot,\cdot\rangle_{V',V}$, and also the inner product in $L^2(D)^{d\times d}$.  Let $f\in V'$ be the forcing function which is assumed to be deterministic. We study the multiscale elasticity equation
\[
-{\partial\over\partial x_j}(a^\ep_{ijkl}(\omega;x)\eps_{kl}(\ue)(\omega;x))=f_i(x),
\]
for all $i=1,\ldots,d$, with the zero boundary condition $\ue=0$ on $\Gamma$ and the traction free condition on $\partial D\setminus\Gamma$. Here $\eps$ denotes the elastic strain tensor
\[
\eps_{ij}(w)=\frac12\left({\partial w_i\over\partial x_j}+{\partial w_j\over\partial x_i}\right)
\]
for functions $w\in H^1(D)^d$. In variational form this problem is written as: Find $\ue\in V$ so that 
\be
\int_Da^\ep(\omega;x)\eps(\ue)(\omega;x):\eps(v)(x)dx=\int_D f(x)\cdot v(x)dx\ \ \forall\,v\in V. 
\label{eq:displacement}
\ee
From Korn's inequality, there is a constant $C$ that only depends on $\alpha$ and $\beta$ in \eqref{eq:boundednesscoerciveness} and the domain $D$ so that
$
\|\ue\|_V\le C\|f\|_{V'}.
$

We study also the Hellinger-Reissner mixed problem that allows for computing the stress tensor $\sigma^\ep(\omega;x)=a^\ep(\omega;x)\eps(\ue)(\omega;x)$ directly. 
 Let $\cH=L^2(D)^{d\times d}_{sym}$. The problem is: Find $(\sigma^\ep(\omega,\cdot),\ue(\omega;\cdot))\in \cH\times V$ so that
\be
		\left\{
	\begin{array}{clrr} %
			((a^\ep)^{-1}(\omega;\cdot)\sigma^\ep(\omega,\cdot),\tau)-(\tau,\epsilon(\ue)(\omega,\cdot)) = 0,\ \  \forall\, \tau \in \cH\\
		  -(\sigma^\ep(\omega,\cdot),\epsilon(v)) = -(f,v),\ \ \forall\, v \in V.
	\end{array}
	\right.
\label{eq:mixed}
\ee
For each $\tau\in \cH$, $a^\ep(\omega;\cdot)\tau$ also belongs to $\cH$. Therefore  the following mixed problem is equivalent to \eqref{eq:mixed}: 
\begin{equation}
		\left\{
	\begin{array}{clrr} %
			(\sigma^\ep(\omega;\cdot),\tau)-(\tau,a^\ep(\omega;\cdot)\epsilon(u^\ep)(\omega;\cdot)) = 0,\ \  \forall\, \tau \in \cH\\
		  -(\sigma^\ep(\omega;\cdot),\epsilon(v)) = -(f,v),\ \ \forall\, v \in V.
	\end{array}
	\right.
\label{eq:mixedeq}
\end{equation}
We study problem \eqref{eq:mixedeq} to use the linear dependence of $a^\ep$ on the random variables $z_m$. 
 From \eqref{eq:boundednesscoerciveness}, the inf-sup conditions for \eqref{eq:mixed} hold uniformly with respect to $\ep$. Therefore 
$
\|\sigma^\ep\|_\cH+\|\ue\|_V\le c\|f\|_{V'}\ \ \forall\,\ep>0.
$
For isotropic materials, the moduli $a$ is written as
\be
a_{ijkl}(\omega;x,{\by})=\mu(\om;x,{\by})(\delta_{ik}\delta_{jl}+\delta_{il}\delta_{jk})+\lambda(\omega;x,{\by})\delta_{ij}\delta_{kl}.
\label{eq:aisotropic}
\ee
We assume a similar structure for the Lam\'e constants $\mu(\omega;x,{\by})$ and $\lambda(\omega;x,{\by})$, i.e.
\begin{equation}\label{eq:mulambda}
\mu(\omega;x,{\by})=\overline \mu(x,{\by})+\sum_{m=1}^\infty z_m(\omega)\mu_m(x,{\by}),
\ \ 
\mbox{and}
\ \ 
\lambda(\omega;x,{\by})=\overline \lambda(x,{\by})+\sum_{m=1}^\infty z_m(\omega)\lambda_m(x,{\by}).
\end{equation}
Let $\gamma_m=\sup_{(x,{\by})\in D\times {\bY}}|\mu_m(x,{\by})|$ and $\delta_m=\sup_{(x,{\by})\in D\times {\bY}}|\lambda_m(x,{\by})|$. We define $\overline\mu_{\min}=\inf_{(x,{\by})\in D\times {\bY}}\overline\mu(x,{\by})$ and $\overline\lambda_{\min}=\inf_{(x,{\by})\in D\times { \bY}}\overline\lambda(x,{\by})$. For the uniform coerciveness and boundedness of the tensor $a_{ijkl}$, we assume:
 There exists a constant $\kappa>0$ such that 
\be
\sum_{m=1}^\infty\gamma_m \leq \frac{\kappa}{\kappa+1}\overline\mu_{\min},\ \ \ 
\sum_{m=1}^\infty\delta_m \leq \frac{\kappa}{\kappa+1}\overline\lambda_{\min},
\mbox{\ \ where\ \ } \overline\mu_{\min}>0, \overline\lambda_{\min}>0.
\label{eq:gammadelta}
\ee

For isotropic materials, the multiplying constant of the best $N$ term convergence rates in the displacement and the Hellinger-Reissner settings depends on the ratio of the Lam\'e constants which is very large when the materials are nearly incompressible. We therefore consider the mixed problem with a penalty term and show that the best $N$ term convergence rate can be established to be independent of this ratio.   
Let $H=L^2(D)$. The mixed problem is:
Find $(u^\ep(w;\cdot),p^\ep(w;\cdot))\in V\times H$ such that
\begin{equation}\label{eq:mixedincomp}
\begin{cases}
2(\mu^\ep(\omega;\cdot)\epsilon(u^\ep(\omega;\cdot)),\epsilon(v(\cdot)))+(\operatorname{div}v(\cdot),p^\ep(\omega;\cdot))=(f,v), \quad \forall v\in V,\\
\ds (\operatorname{div} u^\ep(\omega;\cdot),q(\cdot)) - \left(\frac{1}{\lambda^\ep(\omega;\cdot)}p^\ep(\omega;\cdot),q(\cdot)\right) = 0, \quad \forall q\in H.
\end{cases}
\end{equation}
To employ the linear dependence of $\lambda(\omega;x,{\by})$ on $z_m$, we will consider the equivalent mixed problem
\begin{equation}\label{eq:mixedincompeq}
\begin{cases}
2(\mu^\ep(\omega;\cdot)\epsilon(u^\ep(\omega;\cdot)),\epsilon(v(\cdot)))+(\operatorname{div}v(\cdot),p^\ep(\omega;\cdot))=(f,v), \quad \forall v\in V,\\
\ds \left(\frac{\lambda^\ep(\omega;\cdot)}{\overline\lambda_{\min}}\operatorname{div} u^\ep(\omega;\cdot),q(\cdot)\right)- \left(\frac{1}{\overline\lambda_{\min}}p^\ep(w;\cdot),q(\cdot)\right) = 0, \quad \forall q\in H.
\end{cases}
\end{equation}
We consider the case where ${\rm meas}(\partial D\setminus\Gamma)>0$. Problem \eqref{eq:mixedincomp} has a unique solution (see, e.g. \cite{Braess}) and
$
\|\ue\|_{\bV}+\|p^\ep\|_H+\left\|{p^\ep/(\lambda^\ep)^{1/2}}\right\|_H\le c\|f\|_{V'}
$
where $c$ is independent of $\ep$. 
From \eqref{eq:gammadelta}, $\inf_{x\in D}\lambda^\ep(\omega;x)\ge \bar\lambda_{\min}/(1+\kappa)$ so 
$
\|\ue\|_{\bV}+\|p^\ep\|_{H}\le c\|f\|_{V'}\ \ \forall\,\ep.
$
\begin{remark}\label{rem:incomp} When ${\rm meas}(\partial D\setminus\Gamma)=0$, i.e. $V=H^1_0(D)$, the solution of \eqref{eq:mixedincomp} is bounded with respect to the norm
$
\|\ue\|_\bV+\|p^\ep\|_{H/\IR}+\left\|{p^\ep/(\lambda^\ep)^{1/2}}\right\|_H.
$
This norm is not uniformly equivalent to $\|\ue\|_{\bV}+\|p^\ep\|_H$ when $\bar\lambda_{\min}$ goes to $\infty$. Therefore for nearly incompressible problems, we only consider the case where ${\rm meas}(D\setminus\Gamma)>0$. 
\end{remark}

\section{Deterministic parametric problems} 
To study the law of the solutions of  stochastic problems 
we study parametric problems whose elastic moduli depend on  parameter sequences in $U=[-1,1]^\IN$. We first define the probability space.
\subsection{Probability space}
For the space of parameter sequences $U=[-1,1]^\IN$, we introduce the $\sigma$ algebra $\Sigma_U={\cal B}([-1,1])^\IN$ where ${\cal B}([-1,1])$ is the Borel $\sigma$-algebra on $[-1,1]$. We define the probability measure $\rho$ on $(U,\Sigma_U)$ as 
\[
d\rho(\bz)=\bigotimes_{m=1}^\infty{dz_m\over 2}.
\]
As $dz_m/2$ is a probability measure on $[-1,1]$, $d\rho$ is a probability measure on $U$ so $(U,\Sigma_U,\rho)$ is a probability space. As $z_m$ are independently distributed on $[-1,1]$, for $S=\prod_{m=1}^\infty S_m\subset U$ where $S_m\subset [-1,1]$, 
\[
\rho(S)=\prod_{m=1}^\infty\IP\{\omega:z_m(\omega)\in S_m\}.
\]
\subsection{Parametric deterministic problems}
For $\psi_m(x)$ in \eqref{eq:randommoduli}, we define the deterministic parametric elastic moduli for each $\bz\in U$ as
\[
a(\bz;x,\by)=\bar a(x,\by)+\sum_{m=1}^\infty z_m\psi_m(x,\by).
\]
Conditions \eqref{eq:beta} and \eqref{eq:alpha} guarantee that $a(\bz;x,\by)$ is well defined for all $x\in D$ and $\by\in \bY$ and that
\be
a(\bz;x,\by)\xi:\xi\ge \alpha|\xi|^2,\ \ |a(\bz;x,y)\xi:\zeta|\le \beta|\xi||\zeta|\ \ \forall\,\xi,\zeta\in \IR^{d\times d}_{sym}.
\label{eq:paraboundednesscoerciveness}
\ee
The multiscale parametric elastic moduli are defined as
\[
a^\ep(\bz;x)=a(\bz;x,{x\over\ep_1},\ldots,{x\over\ep_n}).
\] 
We consider the parametric elasticity equation: 
\[
-{\partial\over\partial x_j}(a^\ep_{ijkl}\eps_{kl}(\ue)(\bz;\cdot))=f_i,\ \ i=1,\ldots,d
\]
with the Dirichlet boundary condition $\ue=0$ on $\Gamma$ and the traction free condition on $\partial D\setminus\Gamma$. In variational form, this problem becomes: Find $\ue(\bz;\cdot)\in V$ such that
\be
\int_Da^\ep(\bz;x)\eps(\ue)(\bz;x):\eps(v)(x)dx=\int_D f(x)\cdot v(x)dx\ \forall\,v\in V.
\label{eq:paradisplacement}
\ee
From Lax-Milgram lemma and Korn's inequality, this problem has a unique solution $\ue(\bz;x)$ that satisfies
$
\|\ue(\bz;\cdot)\|_V\le c\|f\|_{V'}
$
where $c$ is independent of $\ep$ and $\bz$. 

We consider the parametric Hellinger-Reissner mixed problem: 
Find $(\sigma^\ep(\bz;\cdot),\ue(\bz;\cdot))\in \cH\times V$ so that
\be
		\left\{
	\begin{array}{clrr} %
			((a^\ep)^{-1}(\bz;\cdot)\sigma^\ep(\bz;\cdot),\tau)-(\tau,\epsilon(\ue)(\bz;\cdot)) = 0,\ \  \forall\, \tau \in H\\
		  -(\sigma^\ep(\bz;\cdot),\epsilon(v)) = -(f,v),\ \ \forall\, v \in V.
	\end{array}
	\right.
\label{eq:paramixed}
\ee
This problem is equivalent to: Find $(\sigma^\ep(\bz;\cdot),\ue(\bz;\cdot))\in \cH\times V$ such that
\begin{equation}
		\left\{
	\begin{array}{clrr} %
			(\sigma^\ep(\bz;\cdot),\tau)-(\tau,a^\ep(\bz;\cdot)\epsilon(\ue)(\bz;\cdot)) = 0,\ \  \forall\, \tau \in \cH\\
		  -(\sigma^\ep(\bz;\cdot),\epsilon(v)) = -(f,v),\ \ \forall\, v \in V.
	\end{array}
	\right.
\label{eq:paramixedeq}
\end{equation}
As the bilinear form $((a^\ep)^{-1}\sigma,\tau)$ on $\cH\times\cH$ is uniformly coercive with respect to $\ep$, the Hellinger-Reissner mixed problem \eqref{eq:paramixed} has a unique solution $(\sigma^\ep,\ue)$ which satisfies
\be
\|\sigma^\ep\|_{\cH}+\|\ue\|_{V}\le c\|f\|_{V'}.
\label{eq:boundsigmaue}
\ee

For nearly incompressible problems, we restrict our consideration to the case where the measure of $\partial D\backslash \Gamma$ is positive. For each $\bz\in U=[-1,1]^\IN$, we consider the parametric Lam\'{e} constants
\[
\mu(\bz;x,{\by})=\overline \mu(x,{\by})+\sum_{m=1}^\infty z_m\mu_m(x,{\by}),\ \ 
\mbox{and}
\ \ 
\lambda(\bz;x,{\by})=\overline \lambda(x,{\by})+\sum_{m=1}^\infty z_m\lambda_m(x,{\by}).
\]
From \eqref{eq:gammadelta}, for all $(x,{\by})\in D\times {\bY}$ and $\bz\in U$ we have
\begin{eqnarray}
\mu_{\min}&:\ds=\frac{1}{1+\kappa}\overline\mu_{\min}\leq\mu(\bz;x,{\by})\leq \sup_{(x,{\by})\in D\times {\bY}}\overline\mu(x,{\by})+\frac{\kappa}{1+\kappa}\overline\mu_{\min}:=\mu_{\max},\label{eq:muminmax}\\
\lambda_{\min}&:\ds=\frac{1}{1+\kappa}\overline\lambda_{\min}\leq \lambda(\bz;x,{\by})\leq \sup_{(x,{\by})\in D\times {\bY}}\overline\lambda(x,{\by})+\frac{\kappa}{1+\kappa}\overline\lambda_{\min}:=\lambda_{\max}.\label{eq:lambdaminmax}
\end{eqnarray}
For problems (\ref{eq:mixedincomp}) and (\ref{eq:mixedincompeq}), we consider the parametric problems: Find $(u^\ep(\bz;\cdot),p^\ep(\bz;\cdot)) \in V\times H$ such that
\begin{equation}\label{eq:paramixedincomp}
\begin{cases}
\ds 2(\mu^\ep(\bz;\cdot)\epsilon(u^\ep(\bz;\cdot)),\epsilon(v(\cdot)))+(\operatorname{div}v(\cdot),p^\ep(\bz;\cdot))=(f,v), \quad \forall v\in V,\\
\ds (\operatorname{div} u^\ep(\bz;\cdot),q(\cdot)) - \left(\frac{1}{\lambda^\ep(\bz;\cdot)}p^\ep(\bz;\cdot),q(\cdot)\right) = 0, \quad \forall q\in L^2(D)
\end{cases}
\end{equation}
and: Find $(u^\ep(\bz;\cdot),p^\ep(\bz;\cdot)) \in V\times H$ such that
\begin{equation}\label{eq:paramixedincompeq}
\begin{cases}
2(\mu^\ep(\bz;\cdot)\epsilon(u^\ep(\bz;\cdot)),\epsilon(v(\cdot)))+(\operatorname{div}v(\cdot),p^\ep(\bz;\cdot))=(f,v), \quad \forall v\in V,\\
\ds\left(\frac{\lambda^\ep(\bz;\cdot)}{\overline\lambda_{\min}}\operatorname{div} u^\ep(\bz;\cdot),q(\cdot)\right)- \left(\frac{1}{\overline\lambda_{\min}}p^\ep(\bz;\cdot),q(\cdot)\right) = 0, \quad \forall q\in L^2(D).
\end{cases}
\end{equation}
Problem \eqref{eq:paramixedincomp} has a unique solution wich satisfies
$
\|\ue\|_V+\|p^\ep\|_H+\left\|{p^\ep/(\lambda^\ep)^{1/2}}\right\|_H\le c\|f\|_{V}.
$
 As $\lambda^\ep\ge \bar\lambda_{\min}/(1+\kappa)$, we deduce that
\be
\|\ue\|_V+\|p^\ep\|_H\le c\|f\|_{V'}.
\label{eq:bounduep}
\ee

To connect the parametric problems
to the stochastic problems
we prove that with respect to the probability measure $(U,\Sigma_U,\rho)$ the solutions are measurable. 
\begin{proposition} \label{prop:measurable}
The solution $\ue(\bz;\cdot)$ of problem \eqref{eq:paradisplacement} as a map from $U$ to $V$ is measurable. The solution $(\sigma^\ep(\bz;\cdot),\ue(\bz;\cdot))$ of problems \eqref{eq:paramixed} and \eqref{eq:paramixedeq} as a map from $U$ to $\cH\times V$ is measurable. The solution $(\ue(\bz;\cdot),p^\ep(\bz;\cdot))$ of problems \eqref{eq:paramixedincomp} and \eqref{eq:paramixedincompeq} as a map from $U$ to $V\times H$ is measurable.  
\end{proposition}
\bproof 
We present the proof for  problem \eqref{eq:paradisplacement}.
The proofs for other problems are similar. 
Let $\bz=(z_1,z_2,\ldots)$ and  $\bz'=(z_1',z_2',\ldots)$ in $U$. 
Let $w=\ue(\bz;\cdot)-\ue(\bz';\cdot)$. We then have
\[
-{\partial\over\partial x_j}(a^\ep_{ijkl}(\bz;\cdot)\eps_{kl}(w))=-{\partial\over\partial x_j}((a^\ep_{ijkl}(\bz';\cdot)-a^\ep_{ijkl}(\bz;\cdot))\eps_{kl}(\ue)(\bz';\cdot)).
\]
From this,
\[
\int_Da^\ep(\bz;x)\eps(w)(x):\eps(w)(x)dx=\int_D(a^\ep(\bz';x)-a^\ep(\bz;x))\eps(\ue)(\bz';x):\eps(w)(x)dx.
\]
As $\ue(\bz;\cdot)$ is uniformly bounded in $V$ with respect to $\bz\in U$, we deduce from \eqref{eq:beta} that 
\[
\|w\|_V\le c\|a^\ep(\bz;\cdot)-a^\ep(\bz';\cdot)\|_{L^\infty(D)^{d^4}}\le c\|a(\bz;\cdot,\cdot)-a(\bz';\cdot,\cdot)\|_{L^\infty(D\times\bY)^{d^4}}.
\]
Thus
\[
\|\ue(\bz;\cdot)-\ue(\bz';\cdot)\|_V\le c\sup_m|z_m-z_m'|\le c\|\bz-\bz'\|_{\ell^\infty(\IN)}.
\]
The mapping $\ue:U\ni\bz\mapsto \ue(\bz;\cdot)\in V$ is thus Lipschitz with respect to the $\ell^\infty(\IN)$ norm, and is therefore measurable. 

\eproof

We therefore have:
\begin{proposition}
For the displacement problem \eqref{eq:displacement}, almost surely, the random solution $\ue(\omega;\cdot)$ can be obtained from the parametric solution $\ue(\bz;\cdot)$ of \eqref{eq:paradisplacement} by
\[
\ue(\omega;\cdot)=\ue(\bz;\cdot)|_{\bz=\bz(\omega)}.
\]
For the stochastic mixed problems \eqref{eq:mixed} and \eqref{eq:mixedeq}, almost surely, 
\[
\sigma^\ep(\omega;\cdot)=\sigma^\ep(\bz;\cdot)|_{\bz=\bz(\omega)},\mbox{\ \ and\ \ }\ue(\omega;\cdot)=\ue(\bz;\cdot)|_{\bz=\bz(\omega)}.
\]
For the stochastic mixed problems \eqref{eq:mixedincomp} and \eqref{eq:mixedincompeq} almost surely
\[
\ue(\omega;\cdot)=\ue(\bz;\cdot)|_{\bz=\bz(\omega)},\mbox{\ \ and\ \ }\ p^\ep(\omega;\cdot)=p^\ep(\bz;\cdot)|_{\bz=\bz(\omega)}.
\]
\end{proposition}
\subsection{Multiscale homegenized problems}
We study the multiscale problems 
via multiscale convergence. We first recall the concept of multiscale convergence developed by Nguetseng \cite{Nguetseng}, Allaire \cite{Allaire} and Allaire and Briane \cite{AllaireBriane}. 
%
\begin{definition}
A sequence $\{w^\ep\}_\ep\subset L^2(D)$ $n+1$-scale converges to a function $w^0\in L^2(D\times\bY)$ if 
\[
\lim_{\ep\to 0}\int_Dw^\ep(x)\phi(x,{x\over\ep_1},\ldots,{x\over\ep_n})dx=\int_D\int_\bY w^0(x,\by)\phi(x,\by)d\by dx\ \ \forall\,\phi\in L^2(D,C_\#(\bY)).
\]
\end{definition}
We have the following results whose proofs can be found in \cite{Nguetseng}, \cite{Allaire} and \cite{AllaireBriane}. 
\begin{proposition}\label{prop:p1}
From a bounded sequence in $L^2(D)$, there exists an $n+1$-scale convergent subsequence.
\end{proposition}
\begin{proposition}\label{prop:p2}
From a bounded sequence $\{w^\ep\}_\ep$ in $H^1(D)$, there exists a subsequence (not renumbered) so that $w^\ep$ converges weakly to $w^0$ in $H^1(D)$, and $n$ functions $w^i\in L^2(D\times Y_1\times\ldots\times Y_{i-1},H^1_\#(Y_i)/\IR)$ ($i=1,\ldots,n$) such that
$
\nabla w^\ep
$
$n+1$-scale converges to $\nabla_xw^0+\nabla_{y_1}w^1+\ldots+\nabla_{y_n}w^n$. 
\end{proposition}
We denote by
$
V_i=L^2(D\times Y_1\times\ldots\times Y_{i-1},H^1_\#(Y_i)/\IR)^d,
$
and
$ 
\bV=V\times V_1\times\ldots\times V_n.
$
The space $\bV$ is equipped with the norm
\[
\|\bv\|_{\bV}=\|v^0\|_{V}+\sum_{i=1}^n\|v^i\|_{V_i}
\]
for $\bv=(v^0,v^1,\ldots,v^n)\in\bV$. For $\bv\in \bV$, for $i,j=1,\ldots,d$, we denote by
\[
\bep_{ij}(\bv)=\frac12\left[\left({\partial v^0_j\over\partial x_i}+{\partial v^1_j\over\partial y_{1i}}+\ldots+{\partial v^n_j\over\partial y_{ni}}\right)+\left({\partial v^0_i\over\partial x_j}+{\partial v^1_i\over\partial y_{1j}}+\ldots+{\partial v^n_i\over\partial y_{nj}}\right)\right].
\]
We have the following Korn type inequality.
\begin{lemma} There is a constant $c$ such that for all $\bv\in \bV$,
\be
\|\bep(\bv)\|_{L^2(D\times\bY)^{d\times d}}\ge c\|\bv\|_\bV.
\label{eq:multiKorn}
\ee
\end{lemma}
\bproof For $\bv=(v^0,v^1,\ldots,v^n)\in\bV$, we have
\[
\|\bep(\bv)\|^2_{L^2(D\times \bY)^{d\times d}}=\|\eps_x(v^0)\|^2_H+\sum_{i=1}^n\|\eps_{y_i}(v^i)\|^2_{L^2(D\times Y_1\times\ldots\times Y_i)^{d\times d}}.
\]
We get the conclusion from Korn's inequality for functions in $V$ and for functions in $(H^1_\#(Y)/\IR)^d$. \eproof

For problem \eqref{eq:paradisplacement}, we have the following result.
\begin{proposition} \label{thm:mshomdisplacement} For each $\bz\in U$, there is $\bu(\bz;\cdot,\cdot)=(u^0,u^1,\ldots,u^n)\in \bV$ such that the parametric solution $\ue(\bz;\cdot)$ of problem \eqref{eq:paradisplacement} converges weakly to $u^0(\bz;\cdot)\in V$ and  $\eps(\ue)$ $n+1$-scale converges to $\bep(\bu)(\bz;\cdot,\cdot)$. The function $\bu$ satisfies the problem
\be
b(\bz;\bu,\bv):=\int_D\int_\bY a(\bz;x,\by)\bep(\bu)(\bz;x,\by):\bep(\bv)(x,\by)d\by dx=\int_Df(x)\cdot v^0(x)dx
\label{eq:msparadisplacement}
\ee
for all $\bv=(v^0,v^1,\ldots,v^n)\in \bV$. Problem \eqref{eq:msparadisplacement} is well-posed. There is a constant $c$ such that
\be
\|\bu(\bz;\cdot,\cdot)\|_\bV\le c\|f\|_{V'}\ \ \forall\,\bz\in U.
\label{eq:boundparabu}
\ee
\end{proposition}
\bproof The proof of this proposition is standard, see, e.g., \cite{Allaire}. As $\ue(\bz;\cdot)$ is uniformly bounded in $V$, 
 there is a subsequence (not renumbered), a function $\bu=(u^0,u^1,\ldots,u^n)\in \bV$ so that $\eps_{ij}(\ue)$ $n+1$-scale converges to $\eps_{xij}(u^0)+\eps_{y_1ij}(u^1)+\ldots+\eps_{y_nij}(u^n)$. By selecting the test function in \eqref{eq:paradisplacement} as 
\[
v(x)=v^0(x)+\sum_{i=1}^n\ep_iv^i(x,{x\over\ep_1},\ldots,{x\over\ep_i})
\]
for $v^0\in{\cal D}(D)^d$ and $v^i(x,y_1,\ldots,y_i)\in {\cal D}(D,C_\#(Y_1\times\ldots\times Y_i))^d$ we get  \eqref{eq:msparadisplacement} from a density argument. The coercivity of $b$ follows from
$
b(\bz;\bu,\bu)\ge c_1\|\bep(\bu)\|_{L^2(D\times\bY)^{d\times d}}^2
$
and \eqref{eq:multiKorn}.
The boundedness of $b$ follows from \eqref{eq:boundednesscoerciveness}. 
 Problem \eqref{eq:msparadisplacement} thus has a unique solution that satisfies \eqref{eq:boundparabu}. \eproof

Let $\bcH=L^2(D\times\bY;\IR^{d\times d}_{sym})$. Let $\bcX=\bcH\times \bV$ with the norm
$
\|(\tau,\bv)\|_{\bcX}=\|\tau\|_\bcH+\|\bv\|_\bV.
$
We define the bilinear forms  $b_1:\bcX\times\bcX\to \IR$ and $b_2:\bcX\times\bcX\to\IR$ as 
\[
b_1((\sigma,\bu),(\tau,\bv))=\int_D\int_\bY \left[a^{-1}(\bz;x,\by)\sigma(\bz;x,\by):\tau(x,\by)-\tau(x,\by):\bep(\bu)(\bz;x,\by)-\sigma(\bz;x,\by):\bep(\bv(x,\by))\right]d\by dx,
\]
\[
b_2((\sigma,\bu),(\tau,\bv))=\int_D\int_\bY \left[\sigma(\bz;x,\by):\tau(x,\by)-a(\bz;x,\by)\tau(x,\by):\bep(\bu)(\bz;x,\by)-\sigma(\bz;x,\by):\bep(\bv(x,\by))\right]d\by dx.
\]
We define the linear form $f:\bV\to \IR$ as
\[
f(\bv)=-\int_Df(x)\cdot v^0(x)dx
\]
for $\bv=(v^0,v^1,\ldots,v^d)\in\bV$. 
We have the following results.
\begin{proposition}\label{thm:mshommixed} For problems \eqref{eq:mixed} and \eqref{eq:mixedeq}, the sequence $\sigma^\ep$ $n+1$- scale converges to $\sigma\in\bcH$ and there is a function $\bu=(u^0,u^1,\ldots,u^n)\in\bV$ such that for all $i,j=1,\ldots,d$, $\eps_{ij}(\ue)$ $n+1$-scale converges to $\bep_{ij}(\bu)$. The functions $\sigma$ and $\bu$ satisfy
\be
b_1(\bz;(\sigma,\bu),(\tau,\bv))=F(\bv),\ \ \forall\ (\tau,\bv)\in\bcX.
\label{eq:b1}
\ee
and
\be 
b_2(\bz;(\sigma,\bu),(\tau,\bv))=F(\bv),\ \ \forall\ (\tau,\bv)\in\bcX.
\label{eq:b2}
\ee
Problems \eqref{eq:b1} and \eqref{eq:b2} possess a unique solution. There are constants $\chi_1$ and $\chi_2$ 
such that
\be
\inf_{(\sigma,\bu)}\sup_{(\tau,\bv)}{b_1(\bz;(\sigma,\bu),(\tau,\bv))\over \|(\sigma,\bu)\|_{\bcX}\|(\tau,\bv)\|_{\bcX}}\ge\chi_1,
\ \ \mbox{and}\ \  
\inf_{(\sigma,\bu)}\sup_{(\tau,\bv)}{b_2(\bz;(\sigma,\bu),(\tau,\bv))\over \|(\sigma,\bu)\|_{\bcX}\|(\tau,\bv)\|_{\bcX}}\ge\chi_2\ \ \forall\,\bz\in U.
\label{eq:infsupb1b2}
\ee

\end{proposition}
\bproof From \eqref{eq:boundsigmaue}, there is  $\sigma(\bz;\cdot,\cdot)\in \bcH$ and  $\bu(\bz;\cdot,\cdot)=(u^0, u^1,\ldots,u^n)\in \bV$ and a susbsequence (not renumbered) so that $\sigma^\ep(\bz;\cdot)$ $n+1$- scale converges to $\sigma$, $\ue\wc u^0$ in $H^1(D)$ and $\eps_{ij}(\ue)$ $n+1$-scale converges to $\eps_{ij}(u^0)+\eps_{y_1ij}(u^1)+\ldots+\eps_{y_nij}(u^n)$. Let $\tau(x,\by)\in C(D\times\bY)^{d\times d}$, $v^0\in {\cal D}(D)^d$ and $v^i\in {\cal D}(D,C_\#(Y_1\times\ldots\times Y_i))^d$. Let $\tau(x,{x\over\ep_1},\ldots,{x\over\ep_n})$ and $v(x)=v^0(x)+\sum_{i=1}^n\ep_iv^i(x,{x\over\ep_1},\ldots,{x\over\ep_i})$ be the test function in \eqref{eq:paramixed}. We then get \eqref{eq:b1}.
 We derive problem \eqref{eq:b2} similarly.

The mapping $\bcH\ni\tau\mapsto a\tau\in \bcH$ is one-to-one so problems \eqref{eq:b1} and \eqref{eq:b2} are equivalent. The inf-sup condition for $b_1$ follows from a standard procedure for mixed elasticity problem (see, e.g., \cite{Braess}), using the uniform coerciveness of $a^{-1}$ and Korn's inequality \eqref{eq:multiKorn}.

\eproof

For a function $\bv=(v^0,v^1,\ldots,v^n)\in\bV$, we define by 
\[
{\bf div}\bv={\rm div}_xv^0+{\rm div}_{y_1}v^1+\ldots+{\rm div}_{y_n}v^n      .
\]
Let $\bH=L^2(D\times\bY)$. For mixed problems \eqref{eq:mixedincomp} and \eqref{eq:mixedincompeq}, we need the following result.
\begin{lemma}\label{lem:divvq}
There is a constant $c_0$ such that  $\forall\,q\in \bH$, there is a function $\bv\in\bV$ so that
\[
{\bf div}\bv(x,\by)=q(x,\by)\ \ 
\mbox{and} 
\ \ 
\|\bv\|_{\bV}\le {\|q\|_{\bH}\over c_0}.
\]
\end{lemma}
\bproof We denote by
\[
Q_0(x)=\int_{\bY}q(x,\by)d\by,\ \ \ \ Q_n(x,y_1,\ldots,y_{n})=q(x,\by)-\int_{Y_n}q(x,\by)dy_n,
\]
\[
Q_i(x,y_1,\ldots,y_{i})=\int_{Y_{i+1}}\ldots\int_{Y_n}q(x,\by)dy_n\ldots dy_{i+1}-\int_{Y_i}\ldots\int_{Y_n}q(x,\by)dy_n\ldots d_{y_i},\ i=1,\ldots,n-1.
\]
Then $q(x,\by)=Q_0(x)+Q_1(x,y_1)+\ldots+Q_n(x,y_1,\ldots,y_n)$. Since ${\rm meas}(\partial D\setminus\Gamma)>0$, there is a function $v^0\in V$ and a constant $c$ such that
$
{\rm div}_x v^0(x)=Q_0(x),\ \ \mbox{and}\ \|v^0\|_{V}\le {\|Q_0\|_H/ c}
$
(see \cite{ErnGuermond} Lemma 4.9 page 181).
For $i=1,\ldots,n$, there is a function $v_i\in V^i$ such that
\[
{\rm div}_{y_i}v^i(x,y_1,\ldots,y_i)=Q_i(x,y_1,\ldots,y_i),\ \ \mbox{and}\ \|v^i(x,y_1,\ldots,y_{i-1},\cdot)\|_{H^1(Y_i)}\le {\|Q_i(x,y_1,\ldots,y_{i-1},\cdot)\|_{L^2(Y_i)}\over c_i}
\]
where the constant $c_i$ is independent of $q$. 
We note that
\[
\|Q_0\|_{L^2(D)}\le \|q\|_{L^2(D\times\bY)},\ \ \mbox{and}\ \|Q_i\|_{L^2(D\times Y_1\times\ldots\times Y_i)}\le 2\|q\|_{L^2(D\times\bY)}.
\]
From these we get the conclusion. \eproof

Let $\bX ={\bV}\times \bH$ with the norm
$
\|({\bu},p)\|_{\bX}= \|{\bu}\|_{\bV}+\|p\|_{\bH}.
$
We define the bilinear forms $b_3, b_4: \bX \times \bX \to \mathbb{R}$ as 
\begin{eqnarray*}
\ds b_3(\bz;({\bu},p),({\bv},q))&=\ds\int_D\int_\bY\left[2\mu(\bz;x,{\by}){\bep}({\bu}(\bz;x,{\by})):{\bep}({\bv}(x,{\by}))+{\bf div}{\bv}(x,{\by})p(\bz;x,{\by})\right]\,d{\by}dx\nonumber\\
 &+\ds\int_D\int_\bY\left[{\bf div}{\bu}(\bz;x,{\by})q(x,{\by})-\frac{1}{\lambda(\bz;x,{\by})}p(\bz;x,{\by})q(x,{\by})\right]\,d{\by}dx
\end{eqnarray*}
and
\begin{eqnarray*}
\ds b_4(\bz;({\bu},p),({\bv},q))&=\ds\int_D\int_\bY\left[2\mu(\bz;x,{\by}){\bep}({\bu}(\bz;x,{\by})):{\bep}({\bv}(x,{\by}))+{\bf div}{\bv}(x,\by)p(\bz;x,{\by})\right]\,d{\by}dx\\
\ds &+\ds\int_D\int_\bY\left[\frac{\lambda(\bz;x,{\by})}{\overline\lambda_{\min}}\operatorname{\bf div}{\bu}(\bz;x,{\by})q(x,{\by})-\frac{1}{\overline\lambda_{\min}}p(\bz;x,{\by})q(x,{\by})\right]\,d{\by}dx.
\end{eqnarray*}

%

 We have the following results. 
\begin{proposition}\label{prop:mshommixedincomp}
For problems (\ref{eq:paramixedincomp}) and (\ref{eq:paramixedincompeq}), the sequence $p^\ep$ $(n+1)$-scale converges to $p\in \bH$ and there is a function ${\bu}=(u^0,u^1,...,u^n)\in {\bV}$ such that $\epsilon_{ij}(u^\ep)$ $(n+1)$-scale converges to ${\bep}({\bu})$. The functions $p$ and ${\bu}$ satisfy:
\begin{equation}\label{eq:b3}
b_3(\bz;({\bu},p),({\bv},q))=-f(\bv), \quad \forall ({\bv},q)\in \bX,
\end{equation}
and:
\begin{equation}\label{eq:b4}
b_4(\bz;({\bu},p),({\bv},q))=-f(\bv), \quad \forall ({\bv},q)\in \bX.
\end{equation}
Problems (\ref{eq:b3}) and (\ref{eq:b4}) are equivalent. There are positive constants $\chi_3$ and $\chi_4$ such that for all $\bz\in U$:
\begin{equation}\label{eq:infsupb3b4}
\inf_{({\bu},p)\in \bX}\sup_{({\bv},q)\in \bX}\frac{b_3(\bz;({\bu},p),({\bv},q))}{\|({\bu},p)\|_{\bX}\|({\bv},q)\|_{\bX}}\geq \chi_3,\ \ 
\mbox{and}\ \ 
\inf_{({\bu},p)\in \bX}\sup_{({\bv},q)\in \bX}\frac{b_4(\bz;({\bu},p),({\bv},q))}{\|({\bu},p)\|_{\bX}\|({\bv},q)\|_{\bX}}\geq \chi_4.
\end{equation}
Problems (\ref{eq:b3}) and (\ref{eq:b4}) have a unique solution. Further, if $\overline\lambda_{\min}>\vartheta$ and $\kappa<\kappa_0$ then $\chi_3$ and $\chi_4$ can be chosen so that they only depend on $\vartheta$ and $\kappa_0$.
\end{proposition}

\bproof
The proof of the limiting problems \eqref{eq:b3} and \eqref{eq:b4} is standard. We use $v(x)=v^0(x)+\sum_{i=1}^n\ep_iv^i(x,\frac{x}{\ep_1},\frac{x}{\ep_2}, ..., \frac{x}{\ep_i})$, where $v^0\in \mathcal D(D)^d$ and $v^i\in \mathcal D(D,C_\#(Y_1\times ... \times Y_i))^d$ and  $q(x,{x\over\ep_1},\ldots,{x\over\ep_n})$ where $q(x,{\by})\in C(D,C_\#(\bY))$ as test functions.

Due to \eqref{eq:lambdaminmax}, $\frac{\lambda}{\overline \lambda_{\min}}$ is uniformly bounded above and below with respect to $x$ and $y$ so for all $q\in \bH$, $\frac{\overline\lambda_{\min}q}{\lambda}\in \bH$. Problems (\ref{eq:b3}) and (\ref{eq:b4}) are equivalent. 
The inf-sup condition of $b_3$ follows the standard procedure for mixed problem with a penality term (see \cite{Braess}) using Lemma \ref{lem:divvq}.  We note that the norm $\|\bv\|_{\bV}+\|q\|_{\bH}+\left\|{q}/{\lambda^\frac{1}{2}}\right\|_{\bH}$ is equivalent to $\|(\bv,q)\|_{\bX}$ uniformly with respect to $\bz\in U$. This norm equivalence is uniform with respect to $\lambda$ if $\overline\lambda_{\min}>\vartheta$ and $\kappa<\kappa_0$ for fixed constants $\vartheta$ and $\kappa_0$. For each $(\bu,p)\in \bX$,
\begin{equation}\label{e1.19}
\sup_{(\bv,p)\in \bX}\frac{b_4(\bz;({\bu},p),({\bv},q))}{\|({\bu},p)\|_{\bX}\|({\bv},q)\|_{\bX}}=\sup_{(\bv,p)\in \bX}\frac{b_3(\bz;({\bu},p),({\bv},\frac{\lambda}{\overline\lambda_{\min}}q))}{\|({\bu},p)\|_{\bX}\|({\bv},\frac{\lambda}{\overline\lambda_{\min}}q)\|_{\bX}}\frac{\|({\bv},\frac{\lambda}{\overline\lambda_{\min}}q)\|_{\bX}}{\|(\bv,q)\|_{\bX}}
\end{equation}
The inf-sup condition for $b_4$ follows from the inequality
$\left\|\frac{\lambda q}{\overline\lambda_{\min}}\right\|_{\bH}\geq \frac{1}{1+\kappa}\|q\|_{\bH}$.
\eproof

We have the following measurablity results.
\begin{proposition} \label{prop:msparameasurable}Solution $u(\bz;\cdot,\cdot)$ of problem \eqref{eq:msparadisplacement} as a map from $U$ to $\bV$ is measurable. Solution $(\sigma(\bz;\cdot,\cdot),u(\bz;\cdot,\cdot))$ of problems \eqref{eq:b1} and \eqref{eq:b2} as a map from $U$ to $\bcX$ is measurable. Solution $(u(\bz;\cdot,\cdot),p(\bz;\cdot,\cdot))$ of problems \eqref{eq:b3} and \eqref{eq:b4} as a map from $U$ to $\bX$ is meassurable. 
\end{proposition}
The proof is similar to that of Proposition \ref{prop:measurable} so we do not present it in details here. 

\section{Approximation of the displacement problem \eqref{eq:msparadisplacement}}
Let $\ubV=L^2(U,\rho;\bV)$. 
From \eqref{eq:boundparabu} and Proposition \ref{prop:measurable}, $\bu$ as a function of $\bz, x$ and $\by$ belongs to $\ubV$. 
We define the bilinear form $B:\ubV\times\ubV\to \IR$ and the linear form $F:\ubV\to \IR$ as 
\[
B(\bu,\bv)=\int_Ub(\bz;\bu,\bv)d\rho(\bz),\ \ 
F(\bv)=\int_U\int_Df(x)\cdot v^0(\bz;x)dxd\rho(\bz),\ \ 
\]
where $b$ is the bilinear form in \eqref{eq:msparadisplacement}.  We consider the problem: Find $\bu\in\ubV$ so that
\be
B(\bu,\bv)=F(\bv)\ \ \ \forall\,\bv\in\ubV.
\label{eq:Bprob}
\ee
As $b$ is uniformly coercive and bounded with respect to $\bz$  so $B$ is bounded and coercive. Problem \eqref{eq:Bprob} has a unique solution. To approximate this solution, 
we  identify a basis of $L^2(U,\rho)$.
\subsection{Orthonormal basis of $L^2(U,\rho)$}
Let $\cF$ be the set of all sequences $\nu=(\nu_j)_{j\ge 1}$ of non-negative integers $\nu_j$ such that only finitely many $\nu_j$ are nonzero. We consider the Legendre polynomials $L_n(t)$ normalized so that
\[
\int_{-1}^1\frac12L_n(t)^2dt=1.
\]
As $L_{\nu_j}(t)=1$ when $j$ is sufficiently large, we can define the multivariate Legendre polynomials as
\[
L_\nu(\bz)=\prod_{j\ge 1}L_{\nu_j}(z_j)
\]
which form an orthonormal basis of $L^2(U,\rho)$. We can therefore write $\bu$ as
\be
\bu=\sum_{\nu\in\cF}\bu_\nu L_\nu(\bz)
\ \ \mbox{where}\ \ 
\bu_\nu=\int_U\bu(\bz;\cdot,\cdot)L_\nu(\bz)d\rho(\bz).
\label{eq:unu}
\ee
\subsection{Semidiscrete Galerkin approximation in $\bz$}
Let $\Lambda$ be a finite subset of $\cF$. We define by
\[
\ubV_\Lambda=\{\bv_\Lambda\in \ubV: \bv_\Lambda(\bz;x,\by)=\sum_{\nu\in\Lambda}\bv_\nu(x,\by)L_\nu(\bz),\bv_\nu\in\bV\}\subset\ubV.
\]
We consider the semidiscrete Galerkin approximation: Find $\bu_\Lambda\in\ubV_\Lambda$ such that
\be
B(\buL,\bv_\Lambda)=F(\bv_\Lambda),\ \ \forall\,\bv_\Lambda\in\ubVL.
\label{eq:GBprob}
\ee
The following error estimate for the semidiscrete Galerkin problem \eqref{eq:GBprob} holds.
\begin{lemma} The solution $\buL$ of problem \eqref{eq:GBprob} satisfies
\[
\|\bu-\buL\|_{\ubV}\le c\left(\sum_{\nu\notin\Lambda}\|\bu_\nu\|_\bV^2\right)^{1/2}.
\]
\end{lemma}
\bproof From Cea's Lemma we have
$
\|\bu-\buL\|_\ubV\le c\inf_{\bv_\Lambda\in \ubV}\|\bu-\bv_\Lambda\|_\bV.
$
Letting $\bv_\Lambda=\sum_{\nu\in \Lambda}\bu_\nu L_\nu$, we get the conclusion using the orthonormality of $L_\nu$.\eproof

\subsection{Bounds for $\|\bu_\nu\|_\bV$}
To get the best $N$ term approximation rate for $\bu$, we deduce bounds for $\|\bu_\nu\|_\bV$. 
For $d_m=\beta_m/\alpha$ where $\beta_m$ and $\alpha$ are the constants in \eqref{eq:beta} and \eqref{eq:boundednesscoerciveness}, we denote by $\boldsymbol d=(d_1,d_2,\ldots)\in\IR^\IN$. We have:
\begin{proposition}\label{prop:boundfordunu}
For the solution $\bu$ of \eqref{eq:msparadisplacement}, there is a constant C independent of $\nu$ such that
\[
\|\bu_\nu(\bz;\cdot,\cdot)\|_\bV\le C|\nu|!{\boldsymbol d}^\nu.
\]
\end{proposition}
\bproof The proof follows the ideas of \cite{CDS1} which has been adapted for stochastic elasticity problems in \cite{XHelasticityrandom}. We show that there is $C_0$ independent of $\nu$ such that 
\be
\|\bep(\partial^\nu_{\bz} \bu)\|_\bcH\le C_0|\nu|!{\boldsymbol d}^\nu\ \ \forall\,\nu\in\cF.
\label{eq:bbepdnu}
\ee 
From \eqref{eq:boundparabu}, we have
$
\|\bep(\bu)(\bz) \|^2_\bcH\le \|f\|_{V'}\|u^0\|_V,
$
so 
$
\|\bep(\bu)(\bz)\|_\bcH\le C_0\ \forall\;\bz\in U.
$
Differentiating both sides of \eqref{eq:msparadisplacement}, we get
\beqas
&&\int_D\int_\bY a(\bz;x,\by)\bep(\partial^\nu_\bz\bu)(\bz;x,\by):\bep(\bv)(x,\by)d\by dx=\\
&&\qquad\qquad\qquad-\sum_{m,\nu_m\ne 0}\nu_m\int_D\int_\bY\psi_m(x,\by)\eps(\partial_\bz^{\nu-e_m}\bu)(\bz;x,\by):\bep(\bv)(x,\by)d\by dx
\eeqas
where $e_m$ is the $m$th unit vector in $\IN^\IN$. 
 Therefore
\[
\|\bep(\partial^\nu_\bz\bu)(\bz;\cdot,\cdot)\|_\bcH\le \sum_m\nu_m{\beta_m\over\alpha}\|\bep(\partial^{\nu-e_m}_\bz \bu)(\bz;\cdot,\cdot)\|_\bcH.
\]
Assuming that \eqref{eq:bbepdnu} holds for $\nu-e_m$, then 
\[
\|\bep(\partial^\nu_\bz\bu)(\bz;\cdot,\cdot)\|_\bcH\le C_0\sum_m\nu_md_m(|\nu|-1)!d^{\nu-e_m}=C_0|\nu|!{\boldsymbol d}^\nu.
\]
Using \eqref{eq:multiKorn}, we get
$
\|\partial^\nu_\bz\bu(\bz;\cdot,\cdot)\|_\bV\le C|\nu|!{\boldsymbol d}^\nu
$
where $C$ is independent of $\nu$.\eproof

Let $\hat d_m=d_m/\sqrt 3$, $\hat {\boldsymbol d}=(\hat d_1,\hat d_2,\ldots)$. We establish the following bounds. 
\begin{proposition}\label{prop:unubound}
For the  expansion \eqref{eq:unu}, for a constant $c$ independent of $\nu$
\[
\|\bu_\nu(\bz;\cdot,\cdot)\|_\bV\le c{|\nu|!\over\nu!}\hat {\boldsymbol d}^\nu.
\]
\end{proposition}
The proof of this proposition uses formula \eqref{eq:unu} and an integration by parts argument following exactly  the procedure in Section 6 of Cohen et al. \cite{CDS1}, using the bound in Proposition \ref{prop:boundfordunu}.
\subsection{Best $N$-term approximation for the solution $\bu$ of problem \eqref{eq:msparadisplacement}}
We deduce the rate of convergence for the best $N$-term approximation for the solution $\bu$. To do so we first  establish the summability of $(\|\bu_\nu\|)_\nu$. We first assume the following.
\begin{assumption}\label{assum:lpsummabilitybetam}
The sequence $(\beta_m)_m$ belongs to $\ell^p(\IN)$ for a constant $0<p<1$.
\end{assumption}

For the summability of $(\|\bu_\nu\|)_\nu$ we employ the following result which is proved in \cite{CDS1}.
\begin{lemma} The sequence $({|\nu|!\over\nu!}\hat {\boldsymbol d}^\nu)_\nu\in \ell^p(\cF)$ if and only if $\|\hat {\boldsymbol d}\|_{\ell^1(\IN)}<1$ and $\hat {\boldsymbol d}\in \ell^p(\IN)$.
\end{lemma}
From \eqref{eq:alpha} and \eqref{eq:ab1} we have that $(1/\alpha)\sum_{m=1}^\infty\beta_m\le\kappa$ so $\sum_{m=1}^\infty \hat d_m<1$ when $\kappa<\sqrt 3$. This, together with Assumption \ref{assum:lpsummabilitybetam}, implies $({|\nu|!\over\nu!}\hat {\boldsymbol d}^\nu)_\nu\in \ell^p(\cF)$. We therefore have
\begin{proposition}
Under Assumption \ref{assum:lpsummabilitybetam}, if $\kappa<\sqrt 3$, then $(\|\bu_\nu\|_\bV)_\nu\in\ell^p(\cF)$. 
\end{proposition} 
 The rate of convergence of the best $N$ term approximation is deduced using the following result.
\begin{lemma} \label{lem:Stechkin} (Stechkin) Let $(b_n)_n$ be a decreasing sequence of positive numbers, then for $0<p<q$ 
\[
\left(\sum_{n>N}b_n^q\right)^{1/q}\le N^{1/q-1/p}\left(\sum_{n\ge 1}b_n^p\right)^{1/p}.
\]
\end{lemma}
With these results, we then have.
\begin{theorem}\label{thm:bestNtermdisplacement} Under Assumption \ref{assum:lpsummabilitybetam} with $\kappa<\sqrt 3$, for any $N$ there is a set $\Lambda\subset\cF$ of cardinality not larger than $N$ such that the error of  the semidiscrete Galerkin approximation problem \eqref{eq:GBprob} satisfies
\[
\|\bu-\bu_\Lambda\|_{\ubV}\le CN^{-s},
\]
where $C$ is independent of $N$ and $s=1/p-1/2$.
\end{theorem}
\bproof Let $\Lambda$ be the set of index sequences $\nu\in \cF$ corresponding to the $N$ terms $\bu_\nu$ in \eqref{eq:unu} with the largest norms $\|\bu_\nu\|_\bV$. We then get the rate of convergence from Lemma \ref{lem:Stechkin}.\eproof

\section{Approximation of mixed problems \eqref{eq:b1} and \eqref{eq:b2}}
We consider the polynomial chaos approximation for problems \eqref{eq:b1} and \eqref{eq:b2} in this section. 
%
%
\subsection{Deterministic parametric Hellinger-Reissner mixed problems}
From Propositions \ref{thm:mshommixed}, the solution of \eqref{eq:b1} and \eqref{eq:b2} satisfies $\sigma\in L^2(U,\rho;\bcH):=\ubcH$ and  $\bu\in L^2(U,\rho;\bV):=\ubV$. Let $\ubcX=\ubcH\times \ubV$. 
We define the bilinear forms $B_1, B_2:\ubcX\times\ubcX\to \IR$ and the linear form $F:\ubcX\to\IR$ as
\beqas
B_1((\sigma,\bu),(\tau,\bv))=\int_Ub_1(\bz;(\sigma,\bu),(\tau,\bv))d\rho(\bz),
\ \ 
\ \ 
B_2((\sigma,\bu),(\tau,\bv))=\int_Ub_2(\bz;(\sigma,\bu),(\tau,\bv))d\rho(\bz),
\eeqas
\[
F((\tau,\bv))=-\int_U\int_Df(x)\cdot v^0(\bz;x)dxd\rho(\bz),
\ \ 
\bv=(v^0,v^1,\ldots,v^n)\in\ubV.
\]
We consider problems: Find $(\sigma,\bu)\in\ubcX$ such that
\be
B_1((\sigma,\bu),(\tau,\bv))=F((\tau,\bv))\ \forall\,(\tau,\bv)\in \ubV
\label{eq:B1prob}
\ee
and 
\be
B_2((\sigma,\bu),(\tau,\bv))=F((\tau,\bv))\ \forall\,(\tau,\bv)\in \ubV. 
\label{eq:B2prob}
\ee
\begin{proposition}\label{prop:Bwellposed}
Problems \eqref{eq:B1prob} and \eqref{eq:B2prob} are equivalent and well-posed.
\end{proposition}
\bproof The proof for the well-posedness of \eqref{eq:B1prob} is standard. As $\tau\mapsto a\tau$ is a one-to-one map from $\ubcH$ to $\ubcH$, these two problems are equivalent. The inf-sup condition for $B_1$ follows from the standard procedure.   
\eproof


\subsection{Semidiscrete Galerkin approximation in $\bz$ of \eqref{eq:B1prob} and \eqref{eq:B2prob}}
We write $\sigma$ and $\bu$ in terms of the multivariate Legendre polynomials $L_\nu$ as
\[
\sigma=\sum_{\nu\in \cF}\sigma_\nu L_\nu,
\mbox{\ and\ } 
\bu=\sum_{\nu\in \cF}\bu_\nu L_\nu,\ \ \sigma_\nu\in\bcH,\ \bu_\nu\in\bV.
\]
Let $\Lambda$ be a subset of $\cF$ of finite cardinality. Let 
\beqas
\ubVL=\{v_\Lambda\in\ubV:\ v_\Lambda=\sum_{\nu\in\Lambda}v_\nu L_\nu,\ v_\nu\in\bV\},\ \ \ 
\ubcHL=\{\tau_\Lambda\in\ubcH:\ \tau_\Lambda=\sum_{\nu\in\Lambda}\tau_\nu L_\nu,\ \tau_\nu\in\bcH\}
\eeqas
and
$
\ubcXL=\ubcHL\times\ubVL.
$
We consider the semidiscrete problems: Find $(\sigma_\Lambda, \buL)\in \ubcXL$ so that
\be
B_1((\sigma_\Lambda,\buL),(\tau_\Lambda,\bv_\Lambda))=F((\tau_\Lambda,\bv_\Lambda))\ \ \forall\,(\tau_\Lambda,\bv_\Lambda)\in\ubcXL
\label{eq:B1Lprob}
\ee
and: Find $(\sigma_\Lambda, \buL)\in \ubcXL$ so that
\be
B_2((\sigma_\Lambda,\buL),(\tau_\Lambda,\bv_\Lambda))=F((\tau_\Lambda,\bv_\Lambda))\ \forall\,(\tau_\Lambda,\bv_\Lambda)\in\ubcXL.
\label{eq:B2Lprob}
\ee
These problems are not equivalent as generally for $\tau\in \ubcHL$, $a\tau$ is not in $\ubcHL$. We therefore establish their well-posedness separately.
\begin{proposition} Problem \eqref{eq:B1Lprob} is well-posed.
\end{proposition}
\bproof The proof of this proposition follows standard proof for saddle point problems. 
\eproof

For problem \eqref{eq:B2Lprob}, we assume that there is a positive constant $\tilde\kappa$ such that 
\be
\sum_{m=1}^\infty\beta_m\le{\tilde\kappa\over 1+\tilde\kappa}{\alpha_0^2\over\alpha_0+\beta_0}
\label{eq:alphastrong}
\ee
 which is stronger than \eqref{eq:alpha}. With this assumption, we have.
\begin{proposition}
If \eqref{eq:alphastrong} holds, problem \eqref{eq:B2Lprob} is well-posed.
\end{proposition}
\bproof We adapt the proof in Xia and Hoang \cite{XHelasticityrandom} for parametric elasticity equations. We define 
\[
\bcH^\eps=\{\zeta\in\bcH:\ \zeta=\bep(\bv),\ \bv\in\bV\}.
\]
From inequality \eqref{eq:multiKorn}, this is a closed subspace of $\bcH$. We define
\[
\ubcH^\eps_\Lambda=\{\zeta_\Lambda\in\ubcH:\ \zeta_\Lambda=\sum_{\nu\in\Lambda}\bep(\bv_\nu)(x,\by)L_\nu(\bz),\ \bv_\nu\in\bV\}.
\]
We define the bilinear forms $\bar{\frak a}_\Lambda:\ubcHL\times\ubcHL\to \IR$ and ${\frak b}_\Lambda:\ubcH_\Lambda\times \ubcH^\ep_\Lambda\to \IR$ as
\[
\bar{\frak a}_\Lambda(\sigma_\Lambda,\tau_\Lambda)=\int_U\int_D\int_\bY\bar a^{-1}(x,\by)\sigma_\Lambda(\bz;x,\by):\tau_\Lambda(\bz;x,\by)d\by dx d\rho(\bz),
\]
\[
{\frak b}_\Lambda(\tau_\Lambda,\zeta_\Lambda)=-\int_U\int_D\int_\bY\tau_\Lambda:\zeta_\Lambda d\by dxd\rho(\bz).
\]
Let $\ubcX_\Lambda^\eps=\ubcH_\Lambda\times\ubcH_\Lambda^\ep$. The bilinear form $\bar B_\Lambda:\ubcXL^\eps\times\ubcXL^\eps$ is defined as 
\[
\bar B_\Lambda((\sigma_\Lambda,\xi_\Lambda),(\tau_\Lambda,\zeta_\Lambda))=\bar {\frak a}_\Lambda(\sigma_\Lambda,\tau_\Lambda)+{\frak b}_\Lambda(\tau_\Lambda,\xi_\Lambda)+{\frak b}_\Lambda(\sigma_\Lambda,\zeta_\Lambda).
\]
 We define the bilinear from $\bar B_{2\Lambda}:\ubcX_\Lambda^\eps\times\ubcX_\Lambda^\ep\to \IR$ as
\beqas
B_{2\Lambda}((\sigma_\Lambda,\xi_\Lambda),(\tau_\Lambda,\zeta_\Lambda))=\bar{\frak a}(\sigma_\Lambda,\tau_\Lambda)+{\frak b}_\Lambda(\tau_\Lambda,\xi_\Lambda)+{\frak b}_\Lambda(\sigma_\Lambda,\zeta_\Lambda)-\\
-\int_U\int_D\int_\bY\tau_\Lambda(\bz;x,\by):\left(\sum_{m=1}^\infty z_m\bar a^{-1}(\bz;x,\by)\psi_m(x,\by)\xi_\Lambda(\bz;x,\by)\right)d\by dxd\rho(\bz).
\eeqas
Let $K$ be the kernel of the map $\ubcHL\ni\tau_\Lambda\mapsto b_\Lambda(\tau_\Lambda,\cdot)\in (\ubcH_\Lambda^\eps)'$. 
From \eqref{eq:abar} we deduce that for all $\bz\in U$, $x\in D$ and $\by\in \bY$ all the eigenvalues of the map $\IR^{d\times d}_{sym}\ni\tau\mapsto\bar a\tau\in\IR^{d\times d}_{sym}$ are not larger than $\beta_0$ so all the eigenvalues of the map $\tau\mapsto\bar a^{-1}\tau$ are not smaller than $1/\beta_0$. Thus for all $\eta\in \IR^{d\times d}_{sym}$, $a^{-1}\eta:\eta\ge \|\eta|^2_{\IR^{d\times d}}/\beta_0$. Then for all $\sigma_\Lambda\in K$
\[
\sup_{\tau_\Lambda\in K}{\bar {\frak a}_\Lambda(\sigma_\Lambda,\tau_\Lambda)\over\|\sigma_\Lambda\|_{\ubcH}\|\tau_\Lambda\|_{\ubcH}}\ge{\bar {\frak a}_\Lambda(\sigma_\Lambda,\sigma_\Lambda)\over\|\sigma_\Lambda\|_{\ubcH}^2} \ge\frac{1}{\beta_0}.
\]
For $\sigma_\Lambda\in \ubcH_\Lambda$, let $\hat\sigma_\Lambda$ be the orthogonal projection of $\sigma_\Lambda$ to $\ubcH_\Lambda^\epsilon$. For all $\xi_\Lambda\in \ubcH_\Lambda^\eps$, we have
\beqas
\sup_{\zeta_\Lambda\in\ubcH_\Lambda^\epsilon}{{\bar B}_\Lambda((\hat\sigma_\Lambda,\xi_\Lambda),(0,\zeta_\Lambda))\over \|\zeta_\Lambda\|_{\ubcH}}=\sup_{\zeta_\Lambda\in\ubcH_\Lambda^\epsilon}{\fb_\Lambda(\hat\sigma_\Lambda,\zeta_\Lambda)\over\|\zeta_\Lambda\|_{\ubcH}}
\ge {\fb_\Lambda(\hat\sigma_\Lambda,\hat\sigma_\Lambda)\over\|\hat\sigma_\Lambda\|_{\ubcH}}=\|\hat\sigma_\Lambda\|_{\ubcH}.
\eeqas
Since $\hat\sigma_\Lambda$ is the orthogonal projection of $\sigma_\Lambda$ to $\ubcH_\Lambda^\ep$, we have $\fb_\Lambda(\sigma_\Lambda,\zeta_\Lambda)=\fb_\Lambda(\hat\sigma_\Lambda,\zeta_\Lambda)$. Therefore
\be
\|\hat\sigma_\Lambda\|_{\ubcH}\le \sup_{(\tau_\Lambda,\zeta_\Lambda)\in\ubcX_\Lambda^\epsilon}{B_{2\Lambda}((\sigma_\Lambda,\xi_\Lambda),(\tau_\Lambda,\zeta_\Lambda))\over\|(\tau_\Lambda,\zeta_\Lambda)\|_{\ubcX_\Lambda^\epsilon}}.
\label{eq:eq}
\ee
Since $\sigma_\Lambda-\hat\sigma_\Lambda\in K$, we have
\beqas
{1\over\beta_0}\|\sigma_\Lambda-\hat\sigma_\Lambda\|_{\ubcH}&\le& \sup_{\tau_\Lambda\in K}{\bar \fa_\Lambda(\sigma_\Lambda-\hat\sigma_\Lambda,\tau_\Lambda)\over\|\tau_\Lambda\|_{\ubcH}}
=\sup_{\tau_\Lambda\in K}{\bar \fa_\Lambda(\sigma_\Lambda-\hat\sigma_\Lambda,\tau_\Lambda)+\fb_\Lambda(\tau_\Lambda,\xi_\Lambda)+\fb_\Lambda(\sigma_\Lambda,0)\over\|(\tau_\Lambda,0)\|_{\ubcX_\Lambda^\epsilon}}\\
&\le&\sup_{(\tau_\Lambda,\zeta_\Lambda)\in\ubcX_\Lambda^\epsilon}{\bar B_\Lambda((\sigma_\Lambda,\xi_\Lambda),(\tau_\Lambda,\zeta_\Lambda))\over \|(\tau_\Lambda,\zeta_\Lambda)\|_{\ubcX_\Lambda^\epsilon}}+\|\bar \fa_\Lambda\|_{\ubcH_\Lambda\times\ubcH_\Lambda\to\IR}\|\hat\sigma_\Lambda\|_{\ubcH}
\eeqas
for all $\xi\in\ubcH_\Lambda^\ep$. Using
\[
\bar B_\Lambda((\sigma_\Lambda,\xi_\Lambda),(\tau_\Lambda,\zeta_\Lambda))\le B_{2\Lambda}((\sigma_\Lambda,\xi_\Lambda),(\tau_\Lambda,\zeta_\Lambda))+{1\over\alpha_0}\left(\sum_{m=1}^\infty\beta_m\right)\|\tau_\Lambda\|_{\ubcH}\|\xi_\Lambda\|_{\ubcH},
\]
we deduce that
\beqas
\|\sigma_\Lambda-\hat\sigma_\Lambda\|_{\ubcH}\le \beta_0\sup_{(\tau_\Lambda,\zeta_\Lambda)\in\ubcX_\Lambda^\epsilon}{B_{2\Lambda}((\sigma_\Lambda,\xi_\Lambda),(\tau_\Lambda,\zeta_\Lambda))\over \|(\tau_\Lambda,\zeta_\Lambda)\|_{\ubcX_\Lambda^\epsilon}}+{\beta_0\over\alpha_0}\left(\sum_{m=1}^\infty\beta_m\right)\|\xi_\Lambda\|_{\ubcH}+{\beta_0\over\alpha_0}\|\hat\sigma_\Lambda\|_{\ubcH}.
\eeqas
Thus together with \eqref{eq:eq} we obtain
\begin{eqnarray}
\|\sigma_\Lambda\|_{\ubcH}&\le& \|\sigma_\Lambda-\hat\sigma_\Lambda\|_{\ubcH}+\|\hat\sigma_\Lambda\|_{\ubcH}\nonumber\\
&\le& (\beta_0+1+{\beta_0\over\alpha_0})\sup_{(\tau_\Lambda,\zeta_\Lambda)\in\ubcX_\Lambda^\epsilon}{B_{2\Lambda}((\sigma_\Lambda,\xi_\Lambda),(\tau_\Lambda,\zeta_\Lambda))\over \|(\tau_\Lambda,\zeta_\Lambda)\|_{\ubcX_\Lambda^\epsilon}}+{\beta_0\over\alpha_0}\left(\sum_{m=1}^\infty\beta_m\right)\|\xi_\Lambda\|_{\ubcH}.
\label{eq:boundsigmaL}
\end{eqnarray}
For $\xi_\Lambda\in\ubcH_\Lambda^\epsilon$ we have
\beqas
\|\xi_\Lambda\|_{\ubcH}&\le& \sup_{\tau_\Lambda\in \ubcH_\Lambda}{\fb_\Lambda(\tau_\Lambda,\xi_\Lambda)\over\|\tau_\Lambda\|_{\ubcH}}\\
&\le&\sup_{\tau_\Lambda\in\ubcH_\Lambda}{\bar \fa_\Lambda(\sigma_\Lambda,\tau_\Lambda)+\fb_\Lambda(\tau_\Lambda,\xi_\Lambda)+\fb_\Lambda(\sigma_\Lambda,0)\over\|(\tau_\Lambda,0)\|_{\ubcX_\Lambda^\epsilon}}+\sup_{\tau_\Lambda\in\ubcH_\Lambda}{\bar \fa_\Lambda(\sigma_\Lambda,\tau_\Lambda)\over\|\tau_\Lambda\|_{\ubcH}}\\
&\le&\sup_{(\tau_\Lambda,\zeta_\Lambda)\in\ubcX_\Lambda^\epsilon}{\bar B_\Lambda((\sigma_\Lambda,\xi_\Lambda),(\tau_\Lambda,\zeta_\Lambda))\over\|(\tau_\Lambda,\zeta_\Lambda)\|_{\ubcX_\Lambda^\epsilon}}+\|\bar \fa_\Lambda\|_{\ubcH_\Lambda\times\ubcH_\Lambda\to{\mathbb R}}\|\sigma_\Lambda\|_{\ubcH}\\
&\le& \sup_{(\tau_\Lambda,\zeta_\Lambda)\in\ubcX_\Lambda^\epsilon}{B_{2\Lambda}((\sigma_\Lambda,\xi_\Lambda),(\tau_\Lambda,\zeta_\Lambda))\over \|(\tau_\Lambda,\zeta_\Lambda)\|_{\ubcX_\Lambda^\epsilon}}+{1\over\alpha_0}\left(\sum_{m=1}^\infty\beta_m\right)\|\xi_\Lambda\|_{\ubcH}+\\
&&{1\over\alpha_0}(\beta_0+1+{\beta_0\over\alpha_0})\sup_{(\tau_\Lambda,\zeta_\Lambda)\in\ubcX_\Lambda^\epsilon}{B_{2\Lambda}((\sigma_\Lambda,\xi_\Lambda),(\tau_\Lambda,\zeta_\Lambda))\over \|(\tau_\Lambda,\zeta_\Lambda)\|_{\ubcX_\Lambda^\epsilon}}+{\beta_0\over\alpha_0^2}\left(\sum_{m=1}^\infty\beta_m\right)\|\xi_\Lambda\|_{\ubcH}\\
&\le&\left(1+{1\over\alpha_0}(\beta_0+1+{\beta_0\over\alpha_0})\right)\sup_{(\tau_\Lambda,\zeta_\Lambda)\in\ubcX_\Lambda^\epsilon}{B_{2\Lambda}((\sigma_\Lambda,\xi_\Lambda),(\tau_\Lambda,\zeta_\Lambda))\over \|(\tau_\Lambda,\zeta_\Lambda)\|_{\ubcX_\Lambda^\epsilon}}+\left({1\over\alpha_0}+{\beta_0\over\alpha_0^2}\right)\left(\sum_{m=1}^\infty\beta_m\right)\|\xi_\Lambda\|_{\ubcH}.
\eeqas
From \eqref{eq:alphastrong}, we deduce that
\be
{1\over 1+\tilde\kappa}\|\xi_\Lambda\|_{\ubcH}\le \left(1+{1\over\alpha_0}(\beta_0+1+{\beta_0\over\alpha_0})\right)\sup_{(\tau_\Lambda,\zeta_\Lambda)\in\ubcX_\Lambda^\epsilon}{B_{2\Lambda}((\sigma_\Lambda,\xi_\Lambda),(\tau_\Lambda,\zeta_\Lambda))\over \|(\tau_\Lambda,\zeta_\Lambda)\|_{\ubcX_\Lambda^\epsilon}}.
\label{eq:boundxiL}
\ee
From \eqref{eq:boundsigmaL} and \eqref{eq:boundxiL}, 
\[
\inf_{(\sigma_\Lambda,\xi_\Lambda)\ubcX_\Lambda^\epsilon}\sup_{(\tau_\Lambda,\zeta_\Lambda)\in\ubcX_\Lambda^\epsilon}{B_{2\Lambda}((\sigma_\Lambda,\xi_\Lambda),(\tau_\Lambda,\zeta_\Lambda))\over \|(\sigma_\Lambda,\xi_\Lambda)\|_{\ubcX_\Lambda^\epsilon}\|(\tau_\Lambda,\zeta_\Lambda)\|_{\ubcX_\Lambda^\epsilon}}\ge c
\]
for a constant $c$ independent of $\Lambda$. 
For $B_2$, we have
\beqas
{B_2((\sigma_\Lambda,\bu_\Lambda),(\tau_\Lambda,\bv_\Lambda))\over \|(\sigma_\Lambda,\bu_\Lambda)\|_{\ubcX_\Lambda}\|(\tau_\Lambda,\bv_\Lambda)\|_{\ubcX_\Lambda}}={B_{2\Lambda}((\sigma_\Lambda,\bep(\bu_\Lambda)),(\bar a\tau_\Lambda,\bep(\bv_\Lambda)))\over \|(\sigma_\Lambda,\bep(\bu_\Lambda))\|_{\ubcX_\Lambda}\|(\bar a\tau_\Lambda,\bep(\bv_\Lambda))\|_{\ubcX_\Lambda}}\cdot {\|(\sigma_\Lambda,\bep(\bu_\Lambda))\|_{\ubcX_\Lambda}\|(\bar a\tau_\Lambda,\bep(\bv_\Lambda))\|_{\ubcX_\Lambda}\over \|(\sigma_\Lambda,\bu_\Lambda)\|_{\ubcX_\Lambda}\|(\tau_\Lambda,\bv_\Lambda)\|_{\ubcX_\Lambda}}.
\eeqas 
From \eqref{eq:multiKorn}, we get the uniform inf-sup condition with respect to $\Lambda$ for bilinear form $B_2$ in \eqref{eq:B2Lprob}.
\eproof

We therefore have the following result on the error estimates.
\begin{proposition}
There is a constant $C$ independent of  $\Lambda$ such that the solution of problems \eqref{eq:B1Lprob} satisfies 
\begin{eqnarray*}
\|(\sigma-\sigma_\Lambda),(\bu-\bu_\Lambda)\|_{\ubcX}\le C\inf_{(\tau_\Lambda,v_\Lambda)\in\ubcX_\Lambda}\|(\sigma-\tau_\Lambda),(\bu-\bv_\Lambda)\|_{\ubcX}
\le C\bigg[\left(\sum_{\nu\notin\Lambda}\|\sigma_\nu\|_{\bcH}^2\right)^{1/2}+\left(\sum_{\nu\notin\Lambda}\|\bu_\nu\|_{\bV}^2\right)^{1/2}\bigg].
\end{eqnarray*}
If condition \eqref{eq:alphastrong} holds, then 
the solution of  problem \eqref{eq:B2Lprob} satisfies 
\begin{eqnarray*}
\|(\sigma-\sigma_\Lambda),(\bu-\bu_\Lambda)\|_{\ubcX}\le C\inf_{(\tau_\Lambda,v_\Lambda)\in\ubcX_\Lambda}\|(\sigma-\tau_\Lambda),(\bu-\bv_\Lambda)\|_{\ubcX}
\le C\bigg[\left(\sum_{\nu\notin\Lambda}\|\sigma_\nu\|_{\bcH}^2\right)^{1/2}+\left(\sum_{\nu\notin\Lambda}\|\bu_\nu\|_{\bV}^2\right)^{1/2}\bigg].
\end{eqnarray*}
\end{proposition}

\subsection{Bounds for $\|\bu_\nu\|_{\bV}$ and $\|\sigma_\nu\|_{\bcH}$}
We deduce in this section explicit bounds for $\|\bu_\nu\|_{\bV}$ and $\|\sigma_\nu\|_{\bcH}$.
 We denote by ${\boldsymbol \delta}=(\delta_1,\delta_2,\ldots)$ where
\[
\delta_m={(1/\alpha_0+\beta_0/\alpha_0^2)\beta_m\over 1-(1/\alpha_0+\beta_0/\alpha_0^2)(\sum_{m=1}^\infty\beta_m)}.
\]
\begin{proposition}\label{prop:bound2}
If \eqref{eq:alphastrong} holds, then there is a constant $C$ independent of $\nu$ such that   
\[
\|\partial_{\bz}^\nu\sigma\|_{\bcH}+\|\partial^\nu_\bz\bu\|_{\bV}\le C|\nu|!{\boldsymbol \delta}^\nu.
\]
\end{proposition}
\bproof
%
For each $\bz\in U$, the solution 
of \eqref{eq:B1prob} and \eqref{eq:B2prob} satisfies
$
\|(\sigma(\bz;\cdot,\cdot), \bu(\bz;\cdot,\cdot))\|_{\bcX}\le c\|f\|_{\bV'}.
$
Let $\bcX^\epsilon=\bcH\times\bcH^\epsilon$ be equipped with the norm
$
\|(\tau,\zeta)\|_{\bcX^\epsilon}=\|\tau\|_{\bcH}+\|\zeta\|_{\bcH}.
$ 
We have
\[
\|(\sigma(\bz;\cdot,\cdot),\bep(\bu)(\bz;\cdot,\cdot))\|_{\bcX^\epsilon}\le c\|f\|_{\bV'}.
\]
Differentiating \eqref{eq:b2}, we have for all $(\tau,\bv)\in \bcX$:
\begin{equation}
\begin{array}{r}
\ds\int_D\int_\bY\partial_\bz^\nu\sigma(\bz;x,\by):\tau(x,\by) d\by dx-\int_D\int_\bY\tau(x,\by):a(\bz;x,\by)\bep(\partial^\nu_\bz \bu)(\bz;x,\by)d\by dx\\
\ds=\sum_{m=1}^\infty\nu_m\int_D\int_\bY\tau(x,\by):\psi_m(x,\by)\bep(\partial_\bz^{\nu-e_m}\bu)(\bz;x,\by)d\by dx
\ds-\int_D\int_\bY\partial^\nu_\bz\sigma(\bz;x,\by):\bep(\bv)(x,\by)d\by dx=0.
\end{array}
\label{eq:deriv}
\end{equation}
As $\bcH\ni\tau\mapsto \bar a(\bz;\cdot,\cdot)\tau\in\bcH$ is one to one, we can  rewrite equations \eqref{eq:deriv} as: $\forall (\tau,\zeta)\in\bcX^\epsilon$
\be
\begin{array}{r}
\ds\int_D\int_\bY\bar a^{-1}\partial^\nu_\bz\sigma(\bz;x,\by):\tau(x,\by)-\int_D\int_\bY\tau(x,\by):\bep(\partial_\bz^\nu\bu)(\bz;x,\by)d\by dx\\
\ds=\sum_{m=1}^\infty\nu_m\int_D\int_\bY\tau(x,\by):\bar a^{-1}(x,\by)\psi_m(x,\by)\bep(\partial_\bz^{\nu-e_m}\bu)(\bz;x,\by)d\by dx\\
\ds\qquad\qquad+\sum_{m=1}^\infty z_m\int_D\int_\bY\tau(x,\by):\bar a^{-1}(x,\by)\psi_m(x,\by)\bep(\partial_\bz^\nu\bu)(\bz;x,\by)d\by dx,\\
\ds-\int_D\int_\bY\partial_\bz^\nu\sigma(\bz;x,\by):\zeta(x,\by)d\by dx=0.
\end{array}
\label{eq:probHHe}
\ee
We have the following inf-sup conditions
\[
\inf_{\sigma\in\bcH}\sup_{\tau\in\bcH}{\displaystyle \int_D\int_\bY\bar a^{-1}(x,\by)\sigma(x,\by):\tau(x,\by)\over\|\sigma\|_{\bcH}\|\tau\|_{\bcH}}\ge{1\over\beta_0},
\mbox{\ and\ } 
\inf_{\zeta\in\bcH^\epsilon}\sup_{\tau\in\bcH}{\displaystyle\int_D\int_\bY\tau(x,\by):\zeta(x,\by)\over\|\tau\|_{\bcH}\|\zeta\|_{\bcH}}\ge 1.
\]
Further for all $\sigma$ and $\tau$ in $\bcH$, 
\[
\left|\int_D\int_\bY\bar a^{-1}(x,\by)\xi(x,\by):\zeta(x,\by)d\by dx\right|\le {1\over\alpha_0}\|\xi\|_{\bcH}\|\zeta\|_{\bcH}.
\] 
Using standard estimates for solutions of saddle point problems (Theorem 2.31 of \cite{ErnGuermond}), we have
\beqas
\|\bep(\partial_\bz^\nu \bu)(\bz;\cdot,\cdot)\|_{\bcH}\le \left(1+{\beta_0\over\alpha_0}\right)\bigg({1\over\alpha_0}\sum_{m=1}^\infty\nu_m\beta_m\|\bep(\partial_\bz^{\nu-e_m}\bu)(\bz;\cdot,\cdot)\|_{\bcH}\qquad\qquad\\
+{1\over\alpha_0}\left(\sum_{m=1}^\infty\beta_m\right)\|\bep(\partial_\bz^\nu \bu)(\bz;\cdot,\cdot)\|_{\bcH}\bigg),
\eeqas
which implies
\beqas
\left(1-\Big({1\over\alpha_0}+{\beta_0\over\alpha_0^2}\Big)\sum_{m=1}^\infty\beta_m\right)\|\bep(\partial_\bz^\nu \bu)(\bz;\cdot,\cdot)\|_{\bcH}\le \left({1\over\alpha_0}+{\beta_0\over\alpha_0^2}\right)\sum_{m=1}^\infty\nu_m\beta_m\|\bep(\partial_\bz^{\nu-e_m}\bu)(\bz;\cdot,\cdot)\|_{\bcH}.
\eeqas
Therefore
\[
\|\bep(\partial_\bz^\nu \bu)(\bz;\cdot,\cdot)\|_{\bcH}\le\sum_{m=1}^\infty\nu_m\delta_m\|\bep(\partial_\bz^{\nu-e_m}\bu)(\bz;\cdot,\cdot)\|_{\bcH}.
\]
By induction, assuming that $\|\bep(\partial_\bz^{\nu-e_m}\bu)(\bz;\cdot,\cdot)\|_{\bcH}\le c(|\nu|-1)!{\boldsymbol \delta}^{\nu-e_m}$ for all $\bz\in U$, then
\[
\|\bep(\partial^\nu_\bz\bu)(\bz;\cdot,\cdot)\|_{\bcH}\le c\sum_m\nu_m(|\nu|-1)!\delta_m{\boldsymbol\delta}^{\nu-e_m}=c|\nu|!{\boldsymbol \delta}^\nu.
\]
From inequality \eqref{eq:multiKorn}, we have
$
\|\partial_\bz^\nu \bu(\bz)\|_{\bV}\le c|\nu|!{\boldsymbol \delta}^\nu.
$
From standard estimates for solutions of saddle point problems and \eqref{eq:probHHe}, we get
\beqas
\|\partial^\nu_\bz\sigma(\bz;\cdot,\cdot)\|_{\bcH}&\le&\beta_0\bigg({1\over\alpha_0}\sum_{m=1}^\infty\nu_m\beta_m\|\bep(\partial_\bz^{\nu-e_m})\bu(\bz;\cdot,\cdot)\|_{\bcH}+
{1\over\alpha_0}\left(\sum_{m=1}^\infty\beta_m\right)\|\bep(\partial_\bz^\nu) \bu(\bz;\cdot,\cdot)\|_{\bcH}\bigg)
\\
&\le& {\beta_0\over\alpha_0}{1-(1/\alpha_0+\beta_0/\alpha_0^2)\sum_{m=1}^\infty\beta_m\over 1/\alpha_0+\beta_0/\alpha_0^2}\sum_{m=1}^\infty\nu_m\delta_mc(|\nu|-1)!{\boldsymbol \delta}^{\nu-e_m}
+{\beta_0\over\alpha_0}\left(\sum_{m=1}^\infty\beta_m\right)c|\nu|!{\boldsymbol \delta}^\nu\\
&\le& C|\nu|!{\boldsymbol \delta}^\nu.
\eeqas

\eproof

Let $\hat \delta_m=\delta_m/\sqrt{3}$ and $\hat{\boldsymbol \delta}=(\hat \delta_1,\hat \delta_2,\ldots)\in{\mathbb R}^{{\mathbb N}}$. We have the estimates:

\begin{proposition}
If condition \eqref{eq:alphastrong} holds, then
\[
\|(\sigma_\nu,\bu_\nu)\|_{\bcX}\le C{|\nu|!\over\nu!}\hat {\boldsymbol \delta}^\nu.
\]
\end{proposition}
The proof of this proposition is similar to that of Proposition \ref{prop:unubound}. 
%
\subsection{Best $N$ term convergence rate}
We deduce the rate of convergence for the semidiscrete Galerkin problems \eqref{eq:B1Lprob} and \eqref{eq:B2Lprob}.
\begin{proposition} \label{prop:psummabilitymixed}
Under Assumption \ref{assum:lpsummabilitybetam}, if \eqref{eq:alphastrong} holds with $\tilde\kappa<\sqrt{3}$ then  
$
(\|(\sigma_\nu,\bu_\nu)\|_{\ubcX})_\nu\in \ell^p(\cF)
$.
\end{proposition}

\bproof
For
 $\left(\displaystyle {\nu!\over|\nu|!}\hat {\boldsymbol \delta}^\nu\right)_{\nu\in{\cal F}}$ to be in $\ell^p(\cF)$ we need $\sum_m\hat \delta_m<1$, i.e. 
\[
\left({1\over\alpha_0}+{\beta_0\over\alpha_0^2}\right)\sum_{m=1}^\infty\beta_m<\sqrt{3}-\sqrt{3}\left({1\over\alpha_0}+{\beta_0\over\alpha_0^2}\right)\sum_{m=1}^\infty\beta_m
\]
so
\[
\sum_{m=1}^\infty\beta_m<{\sqrt{3}\over 1+\sqrt{3}}{\alpha_0^2\over\alpha_0+\beta_0}
\]
which holds when $\tilde\kappa<\sqrt{3}$. 
\eproof

Let $s=1/p-1/2$. 
We then deduce
\begin{theorem}\label{thm:bestNterm12}
Under  Assumption \ref{assum:lpsummabilitybetam}, if   \eqref{eq:alphastrong} holds with  $\tilde\kappa<\sqrt{3}$ then  
for any integer $N$, there is a set $\Lambda\subset\cF$ of cardinality not larger than $N$ such that the error of the semidiscrete Galerkin problems \eqref{eq:B1Lprob} and \eqref{eq:B2Lprob} satisfies
\be
\|\sigma-\sigma_\Lambda\|_{\ubcH}+\|\bu-\bu_\Lambda\|_{\ubV}\le cN^{-s},
\label{eq:bestNterms12}
\ee
where $c$ does not depend on $N$.
\end{theorem}
The proof of this theorem uses Lemma \ref{lem:Stechkin} and Proposition \ref{prop:psummabilitymixed}.
\section{Approximations for mixed problems \eqref{eq:b3} and \eqref{eq:b4}}
We consider approximation for mixed problems for  nearly incompressible materials in this section. We show that the best $N$-term convergence rate does not depend on the ratio of the Lam\'e constants when this ratio goes to $\infty$. As stated in Remark \ref{rem:incomp}, we restrict our consideration to the case where ${\rm meas}(\partial D\setminus\Gamma)>0$. 
\subsection{Deterministic parametric mixed problems for  nearly incompressible materials}
From Proposition \ref{prop:msparameasurable} and the uniform boundedness of $\|(\bu(\bz;\cdot,\cdot), p(\bz;\cdot,\cdot))\|_{\bX}$, the solution of (\ref{eq:b3}) and (\ref{eq:b4}) satisfies $\bu\in L^2(U,\rho;\bV):= \underline {\bV}$ and $p\in L^2(U,\rho;\bH):= \underline {\bH}$. Let $\underline {\bX}=\underline {\bV}\times \underline {\bH}$.
We define the bilinear forms $B_3, B_4: \underline {\bX}\times \underline {\bX}\to \mathbb R$ as 
\[
B_3((\bu,p),(\bv,q))=\int_U b_3(\bz;({\bu},p),({\bv},q))\,d\rho(\bz),
\ \ \ 
B_4((\bu,p),(\bv,q))=\int_U b_4(\bz;({\bu},p),({\bv},q))\,d\rho(\bz).
\]
The linear form $F: \underline {\bX}\to \mathbb R$ is defined as 
$$
F((\bv,q)) = \int_U\int_Df(x)\cdot v^0(\bz;x)\,dxd\rho(\bz).
$$
We consider problems: 
\begin{equation}\label{eq:B3prob}
 \mbox{Find $(\bu, p)\in \underline {\bX}$ such that}\ \  
B_3((\bu,p),(\bv,q))=F((\bv,q)), \quad \forall (\bv,q) \in \underline {\bX},
\end{equation}
and: 
\begin{equation}\label{eq:B4prob}
\mbox{Find $(\bu, p)\in \underline {\bX}$ such that}\ \ 
B_4((\bu,p),(\bv,q))=F((\bv,q)), \quad \forall (\bv,q) \in \underline {\bX}.
\end{equation}
%
\begin{proposition}
Problems \eqref{eq:B3prob} and \eqref{eq:B4prob} are equivalent and possess a unique solution.
\end{proposition}
The proof of this proposition follows standard procedure for saddle point problems with a penalty term. 
\subsection{Semidiscrete Galerkin approximation in $\bz$ of \eqref{eq:B3prob} and \eqref{eq:B4prob}}

As $\bu\in L^2(U,\rho;\bV):= \underline {\bV}$ and $p\in L^2(U,\rho;\bH):= \underline {\bH}$, we can write them as 
\[
\bu=\sum_{\nu\in \mathcal F}\bu_\nu L_\nu, \quad p=\sum_{\nu\in \mathcal F}p_\nu L_\nu,\ \ 
\mbox{where}\ \  \bu_\nu\in \bV,\  p_\nu \in \bH.
\]

For a subset $\Lambda\subset \mathcal F$ of finite cardinality, we denote by 
$$
\underline {\bV}_\Lambda=\left\{\bv_\Lambda \in \bV:~ \bv_\Lambda=\sum_{\nu\in \Lambda}\bv_\nu L_\nu\right\},\ \ \ 
\underline {\bH}_\Lambda=\left\{p_\Lambda \in \bH:~ p_\Lambda=\sum_{\nu\in \Lambda}p_\nu L_\nu\right\},
$$
and $\underline {\bX}_\Lambda=\underline {\bV}_\Lambda\times \underline {\bH}_\Lambda$. We consider the semidiscrete Galerkin approximation for (\ref{eq:B3prob}) and (\ref{eq:B4prob}) as follows: 
\begin{equation}\label{eq:B3Lprob}
\mbox{Find $(\bu_\Lambda, p_\Lambda)\in \underline {\bX}_\Lambda$ so that }\ \ 
B_3((\bu_\Lambda,p_\Lambda),(\bv_\Lambda,q_\Lambda))=F((\bv_\Lambda,q_\Lambda)), \quad \forall (\bv_\Lambda,q_\Lambda) \in \underline {\bX}_\Lambda,
\end{equation}
and: 
\begin{equation}\label{eq:B4Lprob}
\mbox{ Find $(\bu_\Lambda, p_\Lambda)\in \underline {\bX}_\Lambda$ so that }\ \ 
B_4((\bu_\Lambda,p_\Lambda),(\bv_\Lambda,q_\Lambda))=F((\bv_\Lambda,q_\Lambda)), \quad \forall (\bv_\Lambda,q_\Lambda) \in \underline {\bX}_\Lambda.
\end{equation}
We then have the following result.

\begin{proposition}\label{pro:B3Linfsup}
There exists a constant $\chi_3'>0$ independent of $\Lambda$ such that 
\begin{equation}\label{eq:B3Linfsup}
\inf_{(\bu_\Lambda,p_\Lambda)\in \underline {\bX}_\Lambda}\sup_{(\bv_\Lambda,q_\Lambda)\in \underline {\bX}_\Lambda}\frac{B_3((\bu_\Lambda,p_\Lambda),(\bv_\Lambda,q_\Lambda))}{\|(\bu_\Lambda,p_\Lambda)\|_{\underline {\bX}}\|(\bv_\Lambda,q_\Lambda)\|_{\underline {\bX}}}\geq\chi_3'.
\end{equation}
\end{proposition}

The  inf-sup condition (\ref{eq:B3Linfsup}) is the standard result for saddle point problems with a penalty term. For problem (\ref{eq:B4Lprob}) we have the following result.

\begin{proposition}\label{prop:B4Linfsup}
Assume that \eqref{eq:gammadelta} holds with $\kappa\leq \kappa_0$ where $\kappa_0>0$ is a fixed constant. There is a constant $\vartheta_1$ depending on $\kappa_0$, $\mu$ and the domain $D$ such that if $\overline\lambda_{\min}>\vartheta_1$ then 
\begin{equation}\label{eq:B4Linfsup}
\inf_{(\bu_\Lambda,p_\Lambda)\in \underline {\bX}_\Lambda}\sup_{(\bv_\Lambda,q_\Lambda)\in \underline {\bX}_\Lambda}\frac{B_4((\bu_\Lambda,p_\Lambda),(\bv_\Lambda,q_\Lambda))}{\|(\bu_\Lambda,p_\Lambda)\|_{\underline {\bX}}\|(\bv_\Lambda,q_\Lambda)\|_{\underline {\bX}}}\geq\chi_4',
\end{equation}
where $\chi_4'$ is independent of $\bz$. Problem (\ref{eq:B4Lprob}) has a unique solution.

\end{proposition}

\bproof
For $q_\Lambda=\sum_{\nu\in \Lambda}q_\nu L_\nu\in \underline {\bH}_\Lambda$, from Lemma \ref{lem:divvq}, we can choose $\bv_\nu\in \bV$ so that ${\bf div}{\bv}_\nu=q_\nu$ and $\|\bv_\nu\|_{\bV}\leq \frac{1}{c_0}\|q_\nu\|_{\bH}$. Then $\bv_\Lambda=\sum_{\nu\in \Lambda}\bv_\nu L_\nu$ satisfies $\|\bv_\Lambda\|_{\ubV}\leq \frac{1}{c_0} \|q_\Lambda\|_{\ubH}$. From this we deduce
$$
\inf_{q_\Lambda\in\ubHL}\sup_{\bv_\Lambda\in \underline {\bV}_\Lambda}\frac{\int_U\int_D\int_\bY\operatorname{\bf div}{\bv}_\Lambda(\bz;x,{\by})q_\Lambda(\bz;x,{\by})\,d{\by}dxd\rho(\bz)}{\|\bv_\Lambda\|_{\underline {\bV}}\|q_\Lambda\|_{\underline {\bH}}}\geq c_0.
$$
Let $\mu^*$ be a constant such that
\be
\left|\int_U\int_D\int_\bY\mu(\bz;x,{\by}){\bep}({\bu}(\bz;x,{\by})):{\bep}({\bv}(\bz;x,{\by}))\,d\by dxd\rho(\bz)\right|\leq \mu^\ast\|\bu\|_{\ubV}\|\bv\|_{\ubV}\ \  \forall\,\bu,\bv\in\ubV.
\label{eq:infsupdivvq}
\ee
Let $(\bu_\Lambda,p_\Lambda)\in \underline {\bX}_\Lambda$. Adapting the approach in \cite{BraessM2AN}, we first consider the case 
\begin{equation}\label{e1.30}
\|\bu_\Lambda\|_{\ubV} \leq \frac{c_0\|p_\Lambda\|_{\underline {\bH}}}{4\mu^\ast}.
\end{equation}
From \eqref{eq:infsupdivvq}, we get
\beqas
c_0\|p_\Lambda\|_{\underline {\bH}}&\leq& \sup_{\bv_\Lambda\in \underline {\bV}_\Lambda}\frac{\int_U\int_D\int_\bY\operatorname{\bf div}{\bv}_\Lambda(\bz;x,{\by})p_\Lambda(\bz;x,{\by})\,d{\by}dxd\rho(\bz)}{\|\bv_\Lambda\|_{\underline {\bV}}}\\
&=&\sup_{\bv_\Lambda\in \underline {\bV}_\Lambda}\frac{B_4((\bu_\Lambda,p_\Lambda),(\bv_\Lambda,0))}{\|\bv_\Lambda\|_{\underline {\bV}}}-\frac{2\int_U\int_D\int_\bY\mu(\bz;x,{\by}){\bep}({\bu}_\Lambda(\bz;x,{\by})):{\bep}({\bv}_\Lambda(\bz;x,{\by}))\,d\by dxd\rho(\bz)}{\|\bv_\Lambda\|_{\underline {\bV}}}\\
&\leq& \sup_{(\bv_\Lambda,q_\Lambda)\in \underline {\bX}_\Lambda}\frac{B_4((\bu_\Lambda,p_\Lambda),(\bv_\Lambda,q_\Lambda))}{\|(\bv_\Lambda,q_\Lambda)\|_{\underline {\bX}}}+2\mu^\ast\|\bu_\Lambda\|_{\underline {\bV}}\\
&\leq& \sup_{(\bv_\Lambda,q_\Lambda)\in \underline {\bX}_\Lambda}\frac{B_4((\bu_\Lambda,p_\Lambda),(\bv_\Lambda,q_\Lambda))}{\|(\bv_\Lambda,q_\Lambda)\|_{\underline {\bX}}}+\frac{c_0}{2}\|p_\Lambda\|_{\underline {\bV}}.
\eeqas
Using (\ref{e1.30}), we obtain
$$
\sup_{(\bv_\Lambda,q_\Lambda)\in \underline {\bX}_\Lambda}\frac{B_4((\bu_\Lambda,p_\Lambda),(\bv_\Lambda,q_\Lambda))}{\|(\bv,q_\Lambda)\|_{\underline {\bX}}}\geq \frac{c_0}{2}\|p_\Lambda\|_{\underline {\bV}}\geq \min\left\{\frac{c_0}{4},\mu^\ast\right\}\|(\bu_\Lambda,p_\Lambda)\|_{\underline {\bX}}.
$$
We then consider the case 
\begin{equation}\label{e1.31}
\|\bu_\Lambda\|_{\underline {\bV}} > \frac{c_0\|p_\Lambda\|_{\underline {\bH}}}{4\mu^\ast}.
\end{equation}
Let $c_1>0$ be a constant such that
$$
2\int_U\int_D\int_\bY\mu(\bz;x,{\by}){\bep}({\bu}_\Lambda(\bz;x,{\by})):{\bep}({\bu}_\Lambda(\bz;x,{\by}))\,d\by dxd\rho(\bz)\geq c_1\|{\bu}_\Lambda\|_{\underline {\bV}}^2.
$$
Choosing $q_\Lambda=t\operatorname{\bf div}{\bu}_\Lambda$ we have 
\begin{align*}
B_4((\bu_\Lambda,p_\Lambda),(\bu_\Lambda,q_\Lambda))&\geq c_1\|{\bu}_\Lambda\|_{\underline {\bV}}^2-\|p_\Lambda\|_{\underline {\bH}}\|\operatorname{\bf div}{\bu}_\Lambda\|_{\underline {\bH}}+\frac{t}{1+\kappa}\|\operatorname{\bf div}{\bu}_\Lambda\|_{\underline {\bH}}^2-\frac{t}{\overline\lambda_{\min}}\|p_\Lambda\|_{\underline {\bH}}\|\operatorname{\bf div}{\bu}_\Lambda\|_{\underline {\bH}}.
\end{align*}
When $\bar\lambda_{\min}\geq t$ 
$$
B_4((\bu_\Lambda,p_\Lambda),(\bu_\Lambda,q_\Lambda))\geq c_1\|{\bu}_\Lambda\|_{\underline {\bV}}^2+\frac{t}{1+\kappa}\|\operatorname{\bf div}{\bu}_\Lambda\|_{\underline {\bH}}^2-2\|p_\Lambda\|_{\underline {\bH}}\|\operatorname{\bf div}{\bu}_\Lambda\|_{\underline {\bH}}.
$$
Using (\ref{e1.31}), there is $c_3>0$ so that
\begin{align*}
\begin{split}
B_4((\bu_\Lambda,p_\Lambda),(\bu_\Lambda,q_\Lambda))&\geq c_1\|{\bu}_\Lambda\|_{\underline {\bV}}^2+\frac{t}{1+\kappa}\|\operatorname{\bf div}{\bu}_\Lambda\|_{\underline {\bH}}^2-2c_3\|{\bu}_\Lambda\|_{\underline {\bV}}\|\operatorname{\bf div}{\bu}_\Lambda\|_{\underline {\bH}}\\
&\geq c_1\|{\bu}_\Lambda\|_{\underline {\bV}}^2+\frac{t}{1+\kappa}\|\operatorname{\bf div}{\bu}_\Lambda\|_{\underline {\bH}}^2-c_3c_4\|{\bu}_\Lambda\|_{\underline {\bV}}^2-\frac{c_3}{c_4}\|\operatorname{\bf div}{\bu}_\Lambda\|_{\underline {\bH}}^2
\end{split}
\end{align*}
for all positive constants $c_4$. For $c_4=\frac{c_1}{2c_3}$ and $\overline\lambda_{\min}>\vartheta_1=t=\frac{1}{c_1}3c_3^2(1+\kappa)$, we get
$$
B_4((\bu_\Lambda,p_\Lambda),(\bv_\Lambda,q_\Lambda)) \geq c_5\left(\|{\bu}_\Lambda\|_{\underline {\bV}}+\|q_\Lambda\|_{\underline {\bH}}\right)^2
$$
for $c_5=\frac{1}{2}\min\{\frac{c_1}{2},\frac{c_3^2}{c_1t^2}\}$.
Therefore for a constant $c_6$ independent of $\lambda$
$$
\frac{B_4((\bu_\Lambda,p_\Lambda),(\bu_\Lambda,q_\Lambda))}{\|{\bu}_\Lambda\|_{\underline {\bV}}+\|q_\Lambda\|_{\underline {\bH}}}\geq c_5\|{\bu}_\Lambda\|_{\underline {\bV}}\geq c_6(\|{\bu}_\Lambda\|_{\underline {\bV}}+\|p_\Lambda\|_{\underline {\bH}}).
$$
\eproof
\subsection{Bounds for $\bu_\nu$ and $p_\nu$}
We now establish bounds for $\|(u_\nu,p_\nu)\|_\bX$. 

\begin{proposition}\label{pro1.7}
When $\kappa<\kappa_0$, there are positive constants $\vartheta_2, C_i$, $i=1,2,..,5$ depending on $\kappa_0$, $\mu$ and the domain $D$ such that if $\overline\lambda_{\min}>\vartheta_2$, for $0<\zeta<1$, 
\beqas
\|\bep(\partial_\bz^\nu\bu)\|_{ \bcH}+\frac{1}{\lambda^\zeta_{\min}}\|\lambda \operatorname{\bf div}\partial_\bz^\nu\bu\|_{\bH}\leq \left(1+\frac{C_1}{\lambda_{\min}}+\frac{C_4}{\lambda^\zeta_{\min}}\right)\sum_{m=1}^\infty\nu_m\frac{\gamma_m}{\mu_{\min}}\|\bep(\partial^{\nu-e_m}_\bz\bu)\|_{ \bcH}\\
+\left(1+\frac{C_2}{\lambda^{1-\zeta}_{\min}}+\frac{C_3}{\lambda^{2-\zeta}_{\min}}+\frac{C_5}{\lambda_{\min}}\right)\sum_{m=1}^\infty\nu_m\frac{\delta_m}{\lambda_{\min}}\frac{1}{\lambda^\zeta_m}\|\lambda \operatorname{\bf div}\partial_\bz^{\nu-e_m}\bu\|_{\bH}.
\eeqas
\end{proposition}

\bproof
Differentiating \eqref{eq:b4}, we get
\beqas
&&2\int_D\int_\bY\mu(\bz;x,{\by}){\bep}(\partial_\bz^\nu{\bu}(\bz;x,{\by})):{\bep}({\bv}(x,{\by}))\,d\by dx+\int_D\int_\bY\operatorname{\bf div}{\bv}(x,\by)\partial_\bz^\nu p(\bz;x,\by)\,d\by dx\\
&&-\int_D\int_\bY\frac{\lambda(\bz;x,\by)}{\overline\lambda_{\min}}\operatorname{\bf div}\partial_\bz^\nu{\bu}(\bz;x,\by)q(x,\by)\,d\by dx-\int_D\int_\bY\frac{1}{\overline\lambda_{\min}}\partial_\bz^\nu p(\bz;x,\by)q(x,\by)\,d\by dx\\
&&=-2\sum_{m=1}^\infty\nu_m\int_D\int_\bY\mu_m(x,{\by}){\bep}(\partial_\bz^{\nu-e_m}{\bu}(\bz;x,{\by})):{\bep}({\bv}(x,{\by}))\,d\by dx\\
&&-\sum_{m=1}^\infty\nu_m\int_D\int_\bY\frac{\lambda_m}{\overline\lambda_{\min}}\operatorname{\bf div}\partial_\bz^{\nu-e_m}{\bu}(\bz;x,\by)q(x,\by)\,d\by dx
\eeqas
for all $(\bv, q)\in \bX$. 
For each $q\in \bH$,  ${\lambda q}/{\bar\lambda_{\min}}\in \bH$ so we can rewrite this equation as
\begin{align}\label{e1.32}
\begin{split}
&2\int_D\int_\bY\mu(\bz;x,{\by}){\bep}(\partial_\bz^\nu{\bu}(\bz;x,{\by})):{\bep}({\bv}(x,{\by}))\,d\by dx+\int_D\int_\bY\operatorname{\bf div}{\bv}(x,\by)\partial_\bz^\nu p(\bz;x,\by)\,d\by dx\\
&=-2\sum_{m=1}^\infty\nu_m\int_D\int_\bY\mu_m(x,{\by}){\bep}(\partial_\bz^{\nu-e_m}{\bu}(\bz;x,{\by})):{\bep}({\bv}(x,{\by}))\,d\by dx,
\end{split}
\end{align}
\begin{align}\label{e1.33}
\begin{split}
&\int_D\int_\bY\operatorname{\bf div}\partial_\bz^\nu{\bu}(\bz;x,\by)q(x,\by)\,d\by dx\\
&=\int_D\int_\bY\left(\frac{1}{\lambda(\bz;x,\by)}\partial_\bz^\nu p(\bz;x,\by)-\sum_{m=1}^\infty\nu_m\frac{\lambda_m(x,\by)}{\lambda(\bz;x,\by)}\operatorname{\bf div}\partial_\bz^{\nu-e_m}{\bu}(\bz;x,{\by})\right)q(x,\by)\,d\by dx.
\end{split}
\end{align}
We denote by
$$
g(\bz;x,\by)=\frac{1}{\lambda(\bz;x,\by)}\partial_\bz^\nu p(\bz;x,\by)-\sum_{m=1}^\infty\nu_m\frac{\lambda_m(x,\by)}{\lambda(\bz;x,\by)}\operatorname{\bf div}\partial_\bz^{\nu-e_m}{\bu}(\bz;x,{\by})\in \bH.
$$
From Lemma \ref{lem:divvq}, there is a function $\bu_g\in \bV$ such that $\operatorname{\bf div}\bu_g=g$ and 
$
c_0\|\bu_g\|_{\bV}\leq \|g\|_{\bH}.
$
For all $\bv\in \bV$ with $\operatorname{\bf div}\bv=0$ 
\begin{align*}
\begin{split}
&2\int_D\int_\bY\mu(\bz;x,{\by})\Big({\bep}(\partial_\bz^\nu{\bu}(\bz;x,{\by}))-{\bep}({\bu}_g(\bz;x,{\by}))\Big):{\bep}({\bv}(x,{\by}))\,d\by dx\\
&=-2\sum_{m=1}^\infty\nu_m\int_D\int_\bY\mu_m(x,{\by}){\bep}(\partial_\bz^{\nu-e_m}{\bu}(\bz;x,{\by})):{\bep}({\bv}(x,{\by}))\,d\by dx\\
&\hspace{1cm}-2\int_D\int_\bY\mu(\bz;x,{\by}){\bep}({\bu}_g(\bz;x,{\by})):{\bep}({\bv}(x,{\by}))\,d\by dx.
\end{split}
\end{align*}
Letting $\bv=\partial_\bz^\nu{\bu}-\bu_g$ we obtain
$$
\|{\bep}(\partial_\bz^\nu{\bu})-{\bep}({\bu}_g)\|_{ \bcH}\leq \sum_{m=1}^\infty\nu_m\frac{\gamma_m}{\mu_{\min}}\|{\bep}(\partial_\bz^{\nu-e_m}{\bu})\|_{ \bcH}+\frac{\mu_{\max}}{\mu_{\min}}\|{\bep}({\bu}_g)\|_{ \bcH},
$$
which implies
$$
\|{\bep}(\partial_\bz^\nu{\bu})\|_{ \bcH}\leq \sum_{m=1}^\infty\nu_m\frac{\gamma_m}{\mu_{\min}}\|{\bep}(\partial_\bz^{\nu-e_m}{\bu})\|_{ \bcH}+\left(1+\frac{\mu_{\max}}{\mu_{\min}}\right)\|{\bep}({\bu}_g)\|_{ \bcH}.
$$
Let $c_7(d)$ be a constant such that $\|{\bep}({\bv})\|_{ \bcH}\leq c_7\|\bv\|_{\bV}$ for all $\bv\in \bV$. We have 
\begin{eqnarray}
&&\|{\bep}(\partial_\bz^\nu{\bu})\|_{ \bcH}\leq \sum_{m=1}^\infty\nu_m\frac{\gamma_m}{\mu_{\min}}\|{\bep}(\partial_\bz^{\nu-e_m}{\bu})\|_{ \bcH}+\left(1+\frac{\mu_{\max}}{\mu_{\min}}\right)\frac{c_7}{c_0}\|g\|_{\bH}\nonumber\\
&&\ \ \leq \sum_{m=1}^\infty\nu_m\frac{\gamma_m}{\mu_{\min}}\|{\bep}(\partial_\bz^{\nu-e_m}{\bu})\|_{ \bcH}+\left(1+\frac{\mu_{\max}}{\mu_{\min}}\right)\frac{c_7}{c_0\lambda_{\min}}\left(\|\partial^\nu p\|_{\bH}+\sum_{m=1}^\infty\nu_m\delta_m\|\operatorname{\bf div}\partial_\bz^{\nu-e_m}{\bu}\|_{\bH}\right)\nonumber\\
&&\ \ \leq \sum_{m=1}^\infty\nu_m\frac{\gamma_m}{\mu_{\min}}\|{\bep}(\partial_\bz^{\nu-e_m}{\bu})\|_{ \bcH}+\left(1+\frac{\mu_{\max}}{\mu_{\min}}\right)
\frac{c_7}{c_0\lambda_{\min}}\left(\|\partial^\nu p\|_{\bH}+\sum_{m=1}^\infty\nu_m\frac{\delta_m}{\lambda_{\min}}\|\lambda\operatorname{\bf div}\partial_\bz^{\nu-e_m}{\bu}\|_{\bH}\right).\nonumber\\
\label{e1.34}
\end{eqnarray}
From Lemma \ref{lem:divvq}, there is $\bv\in \bV$ such that $\operatorname{\bf div}\bv=\partial_\bz^\nu p$ and $c_0\|\bv\|_{\bV}\leq \|\partial^\nu_\bz p\|_{\bH}$. We deduce from (\ref{e1.32}) 
\begin{align*}
\begin{split}
\|\partial^\nu_\bz p\|_{\bH}^2&\leq 2\mu_{\max}\|{\bep}(\partial_\bz^\nu{\bu})\|_{\bcH}\|{\bep}({\bv})\|_{\bcH}+2\sum_{m=1}^\infty\nu_m\gamma_m\|{\bep}(\partial_\bz^{\nu-e_m}{\bu})\|_{\bcH}\|{\bep}({\bv})\|_{\bcH}\\
&\leq 2\mu_{\max}\frac{c_7}{c_0}\|{\bep}(\partial_\bz^\nu{\bu})\|_{\bcH}\|\partial^\nu_\bz p\|_{\bH}+2\sum_{m=1}^\infty\nu_m\gamma_m\frac{c_7}{c_0}\|{\bep}(\partial_\bz^{\nu-e_m}{\bu})\|_{\bcH}\|\partial^\nu_\bz p\|_{\bH}.
\end{split}
\end{align*}
Thus, using (\ref{e1.34}),
\begin{align}\label{e1.35}
\begin{split}
\|\partial^\nu_\bz p\|_{\bH}&\leq 2\mu_{\max}\frac{c_7}{c_0}\|{\bep}(\partial_\bz^\nu{\bu})\|_{\bcH}+2\sum_{m=1}^\infty\nu_m\gamma_m\frac{c_7}{c_0}\|{\bep}(\partial_\bz^{\nu-e_m}{\bu})\|_{\bcH}\\
&\leq 2\frac{c_7}{c_0}\left(1+\frac{\mu_{\max}}{\mu_{\min}}\right)\sum_{m=1}^\infty\nu_m\gamma_m\|{\bep}(\partial_\bz^{\nu-e_m}{\bu})\|_{ \bcH}+2\frac{c_7^2\mu_{\max}}{c_0^2\lambda_{\min}}\left(1+\frac{\mu_{\max}}{\mu_{\min}}\right)\|\partial^\nu_\bz p\|_{\bH}\\
&\hspace{1cm}+ 2\frac{c_7^2\mu_{\max}}{c_0^2\lambda_{\min}}\left(1+\frac{\mu_{\max}}{\mu_{\min}}\right)\sum_{m=1}^\infty\nu_m\frac{\delta_m}{\lambda_{\min}}\|\lambda\operatorname{\bf div}\partial_\bz^{\nu-e_m}{\bu}\|_{\bH}.
\end{split}
\end{align}
When
$$
\overline\lambda_{\min}>\vartheta_2:=4\mu_{\max}(1+\kappa_0)\left(1+\frac{\mu_{\max}}{\mu_{\min}}\right)\frac{c_7^2}{c_0^2}
$$
we have 
$$
\lambda_{\min}=\frac{\overline\lambda_{\min}}{1+\kappa}>4\mu_{\max}\left(1+\frac{\mu_{\max}}{\mu_{\min}}\right)\frac{c_7^2}{c_0^2}.
$$
Thus 
\begin{align}\label{e1.36}
\begin{split}
\|\partial^\nu_\bz p\|_{\bH}&\leq 4\frac{c_7}{c_0}\left(1+\frac{\mu_{\max}}{\mu_{\min}}\right)\sum_{m=1}^\infty\nu_m\gamma_m\|{\bep}(\partial_\bz^{\nu-e_m}{\bu})\|_{ \bcH}\\
&\hspace{1cm}+4\mu_{\max}\left(1+\frac{\mu_{\max}}{\mu_{\min}}\right)\frac{c_7^2}{c_0^2\lambda_{\min}}\sum_{m=1}^\infty\nu_m\frac{\delta_m}{\lambda_{\min}}\|\lambda\operatorname{\bf div}\partial_\bz^{\nu-e_m}{\bu}\|_{\bH}.
\end{split}
\end{align}
From (\ref{e1.34}) and (\ref{e1.36}), it follows that 
\begin{align}\label{e1.37}
\begin{split}
\|{\bep}(\partial_\bz^\nu{\bu})\|_{ \bcH}&\leq \left(1+4\frac{c_7^2}{c_0^2}\frac{\mu_{\min}}{\lambda_{\min}}\left(1+\frac{\mu_{\max}}{\mu_{\min}}\right)^2\right)\sum_{m=1}^\infty\nu_m{\gamma_m\over\mu_{\min}}\|{\bep}(\partial_\bz^{\nu-e_m}{\bu})\|_{ \bcH}\\
&\hspace{1cm}+\frac{c_7}{c_0\lambda_{\min}}\left(1+\frac{\mu_{\max}}{\mu_{\min}}\right)\left(1+4\frac{c_7^2}{c_0^2}\frac{\mu_{\max}}{\lambda_{\min}}\left(1+\frac{\mu_{\max}}{\mu_{\min}}\right)\right)\cdot \sum_{m=1}^\infty\nu_m\frac{\delta_m}{\lambda_{\min}}\|\lambda\operatorname{\bf div}\partial_\bz^{\nu-e_m}{\bu}\|_{\bH}\\
&=\left(1+\frac{C_1}{\lambda_{\min}}\right)\sum_{m=1}^\infty\nu_m\frac{\gamma_m}{\mu_{\min}}\|{\bep}(\partial_\bz^{\nu-e_m}{\bu})\|_{ \bcH}+\left(\frac{C_2}{\lambda_{\min}}+\frac{C_3}{\lambda_{\min}^2}\right)\sum_{m=1}^\infty\nu_m\frac{\delta_m}{\lambda_{\min}}\|\lambda\operatorname{\bf div}\partial_\bz^{\nu-e_m}{\bu}\|_{\bH};
\end{split}
\end{align}
the constants
$$
C_1=\frac{4c_7^2\mu_{\min}}{c_0^2}\left(1+\frac{\mu_{\max}}{\mu_{\min}}\right)^2, \quad C_2=\frac{c_7}{c_0}\left(1+\frac{\mu_{\max}}{\mu_{\min}}\right), \quad C_3=\frac{4\mu_{\max}c_7^3}{c_0^3}\left(1+\frac{\mu_{\max}}{\mu_{\min}}\right)^2,
$$
do not depend on $\lambda_{\min}$.
From (\ref{e1.35}), we have 
\begin{align}\label{e1.38}
\begin{split}
\|\lambda\operatorname{\bf div}\partial_\bz^{\nu}{\bu}\|_{\bH}&\leq \|\partial_\bz^\nu p\|_{\bH}+\sum_{m=1}^\infty\nu_m\frac{\delta_m}{\lambda_{\min}}\|\lambda\operatorname{\bf div}\partial_\bz^{\nu-e_m}{\bu}\|_{\bH}\\
&\leq 4\frac{c_7\mu_{\min}}{c_0}\left(1+\frac{\mu_{\max}}{\mu_{\min}}\right)\sum_{m=1}^\infty\nu_m\frac{\gamma_m}{\mu_{\min}}\|{\bep}(\partial_\bz^{\nu-e_m}{\bu})\|_{ \bcH}\\
&\hspace{1cm}+\left(1+4\frac{c_7^2}{c_0^2}\frac{\mu_{\max}}{\lambda_{\min}}\left(1+\frac{\mu_{\max}}{\mu_{\min}}\right)\right)\sum_{m=1}^\infty\nu_m\frac{\delta_m}{\lambda_{\min}}\|\lambda\operatorname{\bf div}\partial_\bz^{\nu-e_m}{\bu}\|_{\bH}\\
&=C_4\sum_{m=1}^\infty\nu_m\frac{\gamma_m}{\mu_{\min}}\|{\bep}(\partial_\bz^{\nu-e_m}{\bu})\|_{ \bcH}+\left(1+\frac{C_5}{\lambda_{\min}}\right)\sum_{m=1}^\infty\nu_m\frac{\delta_m}{\lambda_{\min}}\|\lambda\operatorname{\bf div}\partial_\bz^{\nu-e_m}{\bu}\|_{\bH},
\end{split}
\end{align}
where the constants 
$$
C_4=4\frac{c_7\mu_{\min}}{c_0}\left(1+\frac{\mu_{\max}}{\mu_{\min}}\right), \quad C_5=4\frac{\mu_{\max}c_7^2}{c_0^2}\left(1+\frac{\mu_{\max}}{\mu_{\min}}\right)
$$
do not depend on $\lambda_{\min}$. Thus for $0<\zeta<1$  
\begin{align}\label{e1.39}
\begin{split}
&\|{\bep}(\partial_\bz^\nu{\bu})\|_{ \bcH}+\frac{1}{\lambda_{\min}^\zeta}\|\lambda\operatorname{\bf div}\partial_\bz^{\nu}{\bu}\|_{\bH}\le\left(1+\frac{C_1}{\lambda_{\min}}+\frac{C_4}{\lambda_{\min}^\zeta}\right)\sum_{m=1}^\infty\nu_m\frac{\gamma_m}{\mu_{\min}}\|{\bep}(\partial_\bz^{\nu-e_m}{\bu})\|_{ \bcH}\\
&\hspace{1cm}+\left(1+\frac{C_2}{\lambda^{1-\zeta}_{\min}}+\frac{C_3}{\lambda^{2-\zeta}_{\min}}+\frac{C_5}{\lambda_{\min}}\right)\sum_{m=1}^\infty\nu_m\frac{\delta_m}{\lambda_{\min}}\frac{1}{\lambda_{\min}^\zeta}\|\lambda\operatorname{\bf div}\partial_\bz^{\nu-e_m}{\bu}\|_{\bH}.
\end{split}
\end{align}
\eproof

Letting
$$
\widehat d_m=\max\left\{\frac{\gamma_m}{\mu_{\min}}\left(1+\frac{C_1}{\lambda_{\min}}+\frac{C_4}{\lambda_{\min}^\zeta}\right),\frac{\delta_m}{\lambda_{\min}}\left(1+\frac{C_2}{\lambda^{1-\zeta}_{\min}}+\frac{C_3}{\lambda^{2-\zeta}_{\min}}+\frac{C_5}{\lambda_{\min}}\right)\right\}
$$
and $\widehat d=(\widehat d_1, \widehat d_2, ...)$ we have 
\begin{align}\label{e1.40}
\begin{split}
&\|{\bep}(\partial_\bz^\nu{\bu})\|_{ \bcH}+\frac{1}{\lambda_{\min}^\zeta}\|\lambda\operatorname{\bf div}\partial_\bz^{\nu}{\bu}\|_{\bH}\leq \sum_{m=1}^\infty\nu_m\widehat d_m\left(\|{\bep}(\partial_\bz^{\nu-e_m}{\bu})\|_{ \bcH}+\frac{1}{\lambda_{\min}^\zeta}\|\lambda\operatorname{\bf div}\partial_\bz^{\nu-e_m}{\bu}\|_{\bH}\right).
\end{split}
\end{align}
We thus have the following result.

\begin{proposition}\label{pro1.8}
If $\kappa<\kappa_0$ and $\overline\lambda_{\min}>\vartheta_2$ for the constants $\kappa_0$ and $\vartheta_2$ in Proposition \ref{pro1.7}, there is a constant $C_0$ such that the solution $(\bu, p)\in \underline \bX$ of problems (\ref{eq:b3}) and (\ref{eq:b4}) satisfies
$$
\|\partial_\bz^\nu\bu\|_{\bV}+\|\partial^\nu_\bz p\|_{\bH}\leq C_0|\nu|!\hat d^\nu\ \ \forall\,\nu\in\cF.
$$
\end{proposition}

\bproof
From \eqref{eq:infsupb3b4}, there is a constant $c$ independent of $\lambda$ such that 
$\|{\bep}({\bu})\|_{ \bcH}+\|p\|_{\bH}\leq c\|f\|_{\bV'}.$
As $p=\lambda\operatorname{\bf div}{\bu}$ and $\lambda_{\min}=\frac{\overline\lambda_{\min}}{1+\kappa} > \frac{\vartheta_2}{1+\kappa_0}$, there is a constant $C$ independent of $\lambda$ such that
\begin{equation}\label{e1.42}
\|{\bep}({\bu})\|_{ \bcH}+\frac{1}{\lambda_{\min}^\zeta}\|\lambda\operatorname{\bf div}{\bu}\|_{\bH}\leq C.
\end{equation}
From (\ref{e1.40}) we can show by induction that
$
\|{\bep}(\partial_\bz^\nu{\bu})\|_{ \bcH}+\frac{1}{\lambda_{\min}^\zeta}\|\lambda\operatorname{\bf div}\partial_\bz^{\nu}{\bu}\|_{\bH}\leq C_0|\nu|!\widehat d^\nu.
$
The bound for $\|\partial_\bz^\nu{\bu}\|_{\bV}$ follows from \eqref{eq:multiKorn}.
 From \eqref{e1.36}, 
\beqas
\|\partial_\bz^\nu p\|_{\bH}\le 4{c_7\over c_0}\left( 1+{\mu_{\max}\over\mu_{\min}}\right)\sum_{m=1}^\infty\nu_m{\gamma_m\over\mu_{\min}}(|\nu|-1)!{\hat d}^{\nu-e_m}+
4\mu_{\max}\left(1+{\mu_{\max}\over\mu_{\min}}\right){c_7^2\over c_0^2\lambda_{\min}^{1-\zeta}}\sum_{m=1}^\infty\nu_m{\delta_m\over\lambda_{\min}}(|\nu|-1)!{\hat d}^{\nu-e_m}.
\eeqas
We therefore get the bound for $\|\partial_{\bz}^\nu p\|_{\bH}$. 
\eproof

Let  $\hat\fd_m=\frac{\hat d_m}{\sqrt{3}}$ and $\hat\fd=(\hat\fd_1, \hat\fd_2,...)$. We then have the following proposition.

\begin{proposition}\label{pro1.9}
The coefficients $\bu_\nu$ and $p_\nu$ of the gpc expansion for $\bu$ and $p$ satisfy
\begin{equation}\label{e1.43}
\|{\bu}_\nu\|_{\bV}+\|p_\nu\|_{\bH}\leq c\frac{|\nu|!}{\nu!}\hat\fd^\nu.
\end{equation}
\end{proposition}
\subsection{Best $N$ term convergence rate}
To quantify the rate of convergence for the best $N$ term approximation for $\bu$ and $p$, we need to establish the summability property of $\bu_\nu$ and $p_\nu$. We first make the folllowing assumption on the summability of the coefficients of the expantion \eqref{eq:mulambda}.
\begin{assumption}\label{assum:psummabilitygammadelta}
The constants $\gamma_m$ and $\delta_m$ in \eqref{eq:gammadelta} satisfy $\sum_{m=1}^\infty\gamma_m^p<\infty$ and $\sum_{m=1}^\infty\delta_m^p<\infty$.
The sequence $\{\gamma_m\}$ and $\{\delta_m\}$  satisfy 
\be
\sum_{m=1}^\infty\max\left\{{\gamma_m\over\bar\mu_m},{\delta_m\over\bar\lambda_m}\right\}\le{\kappa\over 1+\kappa}.
\label{eq:k1k}
\ee
\end{assumption}

Under this assumption, we have
\begin{proposition}\label{prop:psummabilityunupnu}
Let $\theta>0$. If Assumption \ref{assum:psummabilitygammadelta} holds with $\kappa<{\sqrt{3}\over 1+\theta}$ then  there is a constant $\vartheta_3>0$ depending on $\theta$, $\mu$ and the domain $D$ such that when $\bar\lambda_{\min}>\vartheta_3$, $\{\|\bu_\nu\|_{\bV}\}_\nu$ and $\{\|p_\nu\|_{\bH}\}_\nu$ are in $\ell^p(\mathcal F)$.
\end{proposition}

\bproof
When  $\sum_{m=1}^\infty{\fd}_m<1$, i.e. $\sum_{m=1}^\infty \widehat d_m<\sqrt{3}$
the sequences $(\|u_\nu\|_{\bV})_\nu$ and $(\|p_\nu\|_{\bH})_\nu$ belong to $\ell^p(\mathcal F)$. 
From \eqref{eq:muminmax} and \eqref{eq:lambdaminmax} we have that
\beqas
\sum_{m=1}^\infty {\hat d}_m=\sum_{m=1}^\infty\max\left\{\frac{\gamma_m}{\mu_{\min}}\left(1+\frac{C_1}{\lambda_{\min}}+\frac{C_4}{\lambda_{\min}^\zeta}\right),\frac{\delta_m}{\lambda_{\min}}\left(1+\frac{C_2}{\lambda_{\min}^{1-\zeta}}+\frac{C_3}{\lambda_{\min}^{2-\zeta}}+\frac{C_5}{\lambda_{\min}}\right)\right\}\\
\le\kappa\max\left\{1+\frac{C_1}{\lambda_{\min}}+\frac{C_4}{\lambda_{\min}^\zeta},1+\frac{C_2}{\lambda_{\min}^{1-\zeta}}+\frac{C_3}{\lambda_{\min}^{2-\zeta}}+\frac{C_5}{\lambda_{\min}}\right\}.
\eeqas
If $\bar\lambda_{\min}>\vartheta_3 (\theta)$ for a constant $\theta>0$ then $\lambda_{\min}=\frac{\overline\lambda_{\min}}{1+\kappa}$ satisfies 
\[
\max\left\{\frac{C_1}{\lambda_{\min}}+\frac{C_4}{\lambda_{\min}^\zeta},\frac{C_2}{\lambda_{\min}^{1-\zeta}}+\frac{C_3}{\lambda_{\min}^{2-\zeta}}+\frac{C_5}{\lambda_{\min}}\right\}<\theta.
\]
Thus $\sum_{m=1}^\infty \widehat d_m<\sqrt{3}$ due to $\kappa<\sqrt{3}/(1+\theta)$. From this, we get the conclusion.
\eproof


Let $s=\frac{1}{p}-\frac{1}{2}$. Let $\overline\lambda_{\max}=\sup_{x\in D}\overline\lambda(x)$. We then deduce the best $N$-term convergence rate for the approximations \eqref{eq:B3Lprob} and \eqref{eq:B4Lprob}.

\begin{theorem}\label{thm:bestNtermrateincomp}
Let $\theta>0$. If Assumption \ref{assum:psummabilitygammadelta} holds with $\kappa<\sqrt{3}/(1+\theta)$, then  there is a constant $\vartheta>0$ depending on $\theta, \mu$ such that when $\bar\lambda_{\min}>\vartheta$, for each $N$ there is a set $\Lambda\subset\cF$ with cardinality not more than $N$ and the solution $(\bu, p)$ of problems (\ref{eq:b3}) and (\ref{eq:b4}) and their approximations $(\bu_\Lambda,p_\Lambda)$ in (\ref{eq:B3Lprob}) and (\ref{eq:B4Lprob}) respectively satisfy
\[
\|\bu-\bu_\Lambda\|_{\underline {\bV}}+\|p-p_\Lambda\|_{\underline {\bH}}\leq CN^{-s},
\]
where $C$ depends only on $\frac{\overline\lambda_{\max}}{\overline\lambda_{\min}}$, $\|\{\|\bu_\nu\|_{\bV}\}\|_{\ell^p (\mathcal F)}$ and $\|\{\|p_\nu\|_{\bH}\}\|_{\ell^p(\mathcal F)}$, and in particular it does not depend on the ratio $\lambda_{\max}/\mu_{\min}$ when $\bar\lambda_{\min}\to+\infty$. 
\end{theorem}

\bproof
The approximation $(\buL,p_\Lambda)$ of \eqref{eq:B4Lprob} satisfies 
\begin{align*}
\begin{split}
 \|\bu-\bu_\Lambda\|_{\ubV}+\|p-p_\Lambda\|_{\ubH} \leq \left(1+\frac{\|B_4\|_{\ubX\times \ubX\mapsto \mathbb R}}{\chi_4'}\right)\left(\left\|\sum_{\nu\notin \Lambda}\bu_\nu L_\nu\right\|_{\ubV}+\left\|\sum_{\nu\notin \Lambda}p_\nu L_\nu\right\|_{\ubH}\right),
\end{split}
\end{align*} 
where $\chi_4'$ is the constant in \eqref{eq:B4Linfsup}. The bilinear form $B_4$ satisfies
\[
\left|B_4((\bu_\Lambda,p_\Lambda),(\bv_\Lambda,q_\Lambda))\right|\leq \left(2\mu^\ast+1+\frac{\lambda_{\max}}{\bar\lambda_{\min}}+\frac{1}{\bar\lambda_{\min}}\right)\|(\bu_\Lambda,p_\Lambda)\|_{\ubX}\|(\bv_\Lambda,q_\Lambda)\|_{\ubX}.
\]
Let $\vartheta=\max\{\vartheta_1,\vartheta_2,\vartheta_3\}$ where $\vartheta_1$, $\vartheta_2$ and $\vartheta_3$ are the constants defined above. We have
\begin{align*}
\begin{split}
&\|\bu-\bu_\Lambda\|_{\underline{\bV}}+\|p-p_\Lambda\|_{\underline{\bH}}\\
&\leq \left(1+\frac{1}{\chi_4'}\left(2\mu^\ast+1+\frac{\bar\lambda_{\max}}{\bar\lambda_{\min}}+\frac{\kappa}{1+\kappa}+\frac{1}{\vartheta}\right)\right)\cdot\left(\left\|\sum_{\nu\notin \Lambda}\bu_\nu L_\nu\right\|_{\ubV}+\left\|\sum_{\nu\notin \Lambda}p_\nu L_\nu\right\|_{\ubH}\right)\\
&\leq C\left[\left(\sum_{\nu\notin \Lambda}\|\bu_\nu\|_{\bV}^2\right)^\frac{1}{2}+\left(\sum_{\nu\notin \Lambda}\|p_\nu\|_{\bH}^2\right)^\frac{1}{2}\right],
\end{split}
\end{align*}
where $C$ only depends on $\frac{\bar\lambda_{\max}}{\bar\lambda_{\min}}$.
Letting $\Lambda$ be the set corresponding to the $N$ largest bounds $C|\nu|!\fd^\nu/\nu!$, we get the conclusion. The proof for $B_3$ is similar. 
\eproof

\section{Correctors for the solutions of the multiscale problems}
\subsection{Two-scale problem}
For two scale problems, we can deduce an explicit homogenization rate of convergence. 
For conciseness, we denote $a(\bz;x,\by)$ as $a(\bz;x,y)$.  For  $1\le r,s \le d$, we define the second order symmetric tensor $e^{rs}\in \IR^{d\times d}_{sym}$ as
$e^{rs}_{kl}={1\over 2}(\delta_{rk}\delta_{sl}+\delta_{rl}\delta_{sk})$.
The homogenized elastic moduli is determined by
\be
a^0_{ijkl}=\int_Y a(\bz;x,y)(e^{ij}+\eps_y(N^{ij})):(e^{kl}+\eps_y(N^{kl}))dy.
\label{eq:a0}
\ee
where $N^{rs}$ is the solution of the cell problem
\be
\int_Y a(\bz;x,y)(e^{rs}+\eps_y(N^{rs})):\eps_y(\phi)dy=0,\ \ \forall\,\phi\in H^1_\#(Y)^d.
\label{eq:2scell}
\ee
We have the following homogenization rate of convergence.
\begin{proposition}\label{prop:homrateconv}
If $u^0\in L^\infty(U;H^2(D)^d)$, $N^{rs}\in L^\infty(U;C^1(\bar D,C^1(\bar Y)\cap H^2(Y)))^d$ and $\partial D$ is Lipschitz, then
\be
\left\|\ue(\bz;\cdot)-u^0(\bz;\cdot)-\ep u^1\left(\bz;\cdot,{\cdot\over\ep}\right)\right\|_{L^\infty(U;H^1(D)^d)}\le c\ep^{1/2}.
\label{eq:homerror}
\ee
\end{proposition}
The proof of this proposition is similar to that for the non-parametric case which can be found in \cite{XHelasticity}. The uniform constant $c$ in the homogenization rate with respect to the parameters is due to the uniform boundednes and coerciveness of the elasticity moduli $a$, and the uniform regularity of $u^0$ and $N^{rs}$. 
To have the required regularity for $N^{rs}$, we make the following assumption..
\begin{assumption}\label{assum:psim2s} The fourth order tensors $\bar a$ and $\psi_m$ in \eqref{eq:randommoduli} belongs to $C^2(D,C^2_\#(Y))^{d^4}$ such that
\[
\sum_{m=1}^\infty\|\psi_m\|_{C^2(D,C^2_\#(Y))^{d^4}}<\infty.
\]
\end{assumption} 
The elastic tensor $a(\bz;\cdot,\cdot)$ is then  uniformly bounded in $C^2(D,C^2_\#(Y))^{d^4}$. We then have:
\begin{lemma}
 \label{lem:Nrsregularity}
Under Assumption \ref{assum:psim2s}, $N^{rs}\in L^\infty(U;C^1(\bar D, C^1(\bar Y)))^d$.
\end{lemma}

For the homogenization error estimate \eqref{eq:homerror} to hold, we need the following result.
\begin{lemma}\label{lem:u0LinftyH2}
Assume that $\partial D$ belongs to the $C^1$ class and $f\in L^2(D)^d$. Then $u^0\in L^\infty(U;H^2(D)^d)$.
\end{lemma}
The proofs of Lemmas \ref{lem:Nrsregularity} and \ref{lem:u0LinftyH2} use eliptic regularity (Theorems 4.16 and 4.18 of \cite{Mclean}). We refer to \cite{XHelasticity} for details. 
From this we deduce.
\begin{proposition}\label{prop:prop} Assume that $\partial D$ belongs to the $C^1$ class and $f\in L^2(D)^d$. Under Assumption \ref{assum:psim2s}, 
\[
\left\|\nabla\ue(\bz;\cdot)-\left[\nabla u^0(\bz;\cdot)+\nabla_yu^1\left(\bz;\cdot,{\cdot\over\ep}\right)\right]\right\|_{L^2(U,\rho;L^2(D)^{d\times d})}\le c\ep^{1/2}.
\]
\end{proposition}
To deduce an approximation for the solution $\ue$ of the multiscale parametric problem \eqref{eq:paradisplacement}, we introduce the  operator ${\cal U}^\ep: L^1(D\times Y)\to L^1(D)$ which is defined as 
\[
{\cal U}^\ep(\Phi)(x)=\int_Y\Phi\left(\ep\left[{x\over\ep}\right]+\ep t,\left\{{x\over\ep}\right\}\right)dt
\]
where $[\cdot]$ denotes the integer part with respect to $Y$ and $\{\cdot\}=\cdot-[\cdot]$. Let $D^\ep$ be a $2\ep$ neighbourhood of $D$, we have:
\begin{lemma}\label{lem:l1}
For $\Phi\in L^1(D\times Y)$,
\[
\int_{D^\ep}{\cal U}^\ep(\Phi)(x)dx=\int_D\int_Y\Phi(x,y)dydx.
\]
\end{lemma} 
A proof of this Lemma can be found in \cite{CDG}.
We then have the following result.
\begin{lemma}\label{lem:l2}
If $u^0\in L^\infty(U;H^2(D))^{d}$ and $N^{rs}\in L^\infty(U;C^1(D,C^1(Y)))^d$, then 
\[
\sup_{\bz\in U}\int_D\left|\nabla_yu_1(\bz;x,{x\over\ep})-{\cal U}^\ep(\nabla_yu_1(\bz;\cdot,\cdot))(x)\right|^2dx\le c\ep^2
\]
where the constant $c$ is independent of $\ep$ and $\bz\in U$ .
\end{lemma}
\bproof
This Lemma is essentially Lemma 5.5 in Hoang and Schwab \cite{HSmultirandom}. It relies on the fact that 
\[
\int_D\int_Y\left|\eps(u^0)(\bz;x)-\eps(u^0)\left(\bz;\ep\left[{x\over\ep}\right]+\ep t \right)\right|^2dtdx\le c\ep^2
\]
as $\eps(u^0)\in H^1(D)^{d\times d}$; $c$ only depends on $\|\eps(u^0)\|_{H^1(D)^{d\times d}}$. The proof for this is quite technical so we refer to \cite{HSmultirandom} for details. Further as $N^{rs}\in L^\infty(U;C^1(D,C^1(Y)))^d$,
\[
{\rm ess}\sup_{\bz\in U}\sup_{t\in Y}\left|\nabla_yN^{rs}\left(\bz;x,{x\over\ep}\right)-\nabla_yN^{rs}\left(\bz;\ep\left[{x\over\ep}\right]+\ep t,{x\over\ep}\right)\right|\le c\ep.
\]
From these we get the conclusion.\eproof

For the mixed problem \eqref{eq:paramixed}, we have the following approximations.
\begin{theorem}
Assume that the boundary $\partial D$ belongs to the class $C^1$, and $f\in L^2(D)^d$.  If condition \eqref{eq:alphastrong} holds, then for each $ N$ there is a set $\Lambda_N\subset\cF$ of cardinality not more than $N$ such that the solution of the approximating problems \eqref{eq:B1Lprob} and \eqref{eq:B2Lprob} satisfy
\[
\|\nabla\ue-[\nabla u^0_{\Lambda_N}+{\cal U}^\ep(\nabla_yu^1_{\Lambda_N})]\|_{L^2(U,\rho;L^2(D)^{d\times d})}+\|\sigma^\ep-{\cal U}^\ep(\sigma_{\Lambda_N})\|_{L^2(U,\rho;L^2(D)^{d\times d})}\le c(\ep^{-1/2}+N^{-s}).
\]
\end{theorem}
\bproof From Proposition \ref{prop:prop}, we have
\[
\|\sigma^\ep-a^\ep[\eps(u^0)+\eps_y(u^1)(\cdot;\cdot,{\cdot\over\ep})]\|_{L^2(U,\rho;L^2(D)^{d\times d})}\le c\ep^{1/2}
\ \ \ 
\mbox{i.e.}\ \ 
\|\sigma^\ep-\sigma(\cdot;\cdot,{\cdot\over\ep})\|_{L^2(U,\rho;L^2(D)^{d\times d})}\le c\ep^{1/2}.
\]
From Lemma \ref{lem:l2}, we have
$
\sup_{\bz\in U}\int_D|\eps_y(u^1)(\bz;x,{x\over\ep})-{\cal U}^\ep(\eps_y(u^1)(\bz;\cdot,\cdot))(x)|^2dx\le c\ep^2.
$
From the proof of Lemma \ref{lem:l2}, as $\eps(u^0)(\bz;\cdot)$ is uniformly bounded in $H^1(D)^{d\times d}$ 
\[
\sup_{\bz\in U}\int_D\left|\eps(u^0)(\bz;x)-{\cal U}^\ep(\eps(u^0))(\bz;x)\right|^2dx\le c\ep^2.
\]
These imply
$
\sup_{\bz\in U}\int_D\left|\eps(u^0)(\bz;x)+\eps_y(u^1)(\bz;x,{x\over\ep})-{\cal U}^\ep(\eps(u^0)(\bz;\cdot)+\eps_y(u^1)(\bz;\cdot,\cdot))(x)\right|^2dx\le c\ep^2.
$
Therefore
\[
\sup_{\bz\in U}\int_D\left|\sigma(\bz;x,{x\over\ep})-a(\bz;x,{x\over\ep}){\cal U}^\ep(\eps(u^0)(\bz;\cdot)+\eps_y (u^1)(\bz;\cdot,\cdot))(x)\right|^2dx\le c\ep^2.
\]
Since $a\in L^\infty(U;C^1(D,C^1(Y)))$,
$
\sup_{\bz\in U}\sup_{x\in D}|a(\bz;x,{x\over\ep})-{\cal U}^\ep(a)(\bz;x)|\le c\ep.
$
 Using ${\cal U}^\ep(a){\cal U}^\ep(\eps(u^0)+\eps_y(u^1))={\cal U}^\ep(a(\eps(u^0)+\eps_y(u^1)))={\cal U}^\ep(\sigma)$, we deduce that
\[
\sup_{\bz\in U}\int_D|\sigma(\bz;x,{x\over\ep})-{\cal U}^\ep(\sigma)(x)|^2dx\le c\ep^2.
\]
From Theorem \ref{thm:bestNterm12}, by choosing $\Lambda_N$ as the set corresponding to the indices $\nu$ with the largest $\|\bu_\nu\|_{\bV}+\|\sigma_\nu\|_{\bcH}$ we have
\beqas
\|{\cal U}^\ep(\nabla_yu^1-\nabla_yu^1_{\Lambda_N})\|_{L^2(U,\rho;L^2(D)^{d\times d})}+\|{\cal U}^\ep(\sigma-\sigma_{\Lambda_N})\|_{L^2(U,\rho;L^2(D)^{d\times d})}\le \\
\|\nabla_yu^1-\nabla_yu^1_{\Lambda_N}\|_{L^2(U,\rho;L^2(D\times Y)^{d\times d})}+\|\sigma-\sigma_{\Lambda_N}\|_{L^2(U,\rho;L^2(D\times Y)^{d\times d})}\le cN^{-s}.
\eeqas
Therefore 
\[
\|\nabla u^\ep-[\nabla u^0_{\Lambda_N}+{\cal U}^\ep_n(\nabla_yu^1_{\Lambda_N})]\|_{L^2(U,\rho;L^2(D)^{d\times d})}+\|\sigma^\ep-{\cal U}^\ep(\sigma_{\Lambda_N})\|_{L^2(U,\rho;L^2(D)^{d\times d})}\le c(\ep^{1/2}+N^{-s}).
\]
\eproof
\subsection{Two-scale nearly incompressible problem}
The constant $c$ in estimate \eqref{eq:homerror} depends explicitly on the ratio ${\sup_{x,\by}\lambda/\inf_{x,\by}\mu}$ which is very large when the material is nearly incompressible. In this section, we deduce a homogenization error rate that does not depend explicitly on this ratio. 
 Let 
\[
||\alpha||=\max_{ijklr}\|\alpha_{jr}^{ikl}\|_{L^\infty(U;C^1(\bar D,C(\bar Y)))},\ \ ||N||=\max_{r,s}\|N^{rs}\|_{L^\infty(U;C^1(\bar D,C^1(\bar Y)\cap H^2(Y)))},\ \ \mbox{and}\ ||u_0||=\|u_0\|_{H^2(D)}
\]
where $\alpha$ is the tensor defined in \eqref{eq:alphatensor} below. We then have: 
\begin{proposition}
If $u^0\in L^\infty(U;H^2(D)^d)$ and $N^{rs}\in L^\infty(U;C^1(\bar D,C^1(\bar Y)\cap H^2(Y)))^d$, then there are constants $c_1=c_1(||\alpha||,||N||,||u^0||)$ and $c_2=c_2(||N||,||u^0||)$ such that
\be
\left\|\ue(\bz;\cdot)-u^0(\bz;\cdot)-\ep u^1(\bz;\cdot,{\cdot\over\ep})\right\|_{H^1(D)^d}\le c_1\ep+c_2\ep^{1/2}.
\label{eq:homerrornearlyincomp}
\ee
\end{proposition}
\bproof
We consider the cell problem
\begin{equation}
		\left\{
	\begin{array}{clrr} %
\ds\int_Y\left[2\mu(\bz;x,y)(e^{rs}+\eps_y(N^{rs}(\bz;x,y))):\eps(\phi)+{\rm div}_y\phi(y)p^{rs}(\bz;x,y)\right]dy=0\\
\ds\int_Y\left[(e^{rs}_{ii}+{\rm div}_yN^{rs}(\bz;x,y))q(x,y)-{1\over\lambda(\bz;x,y)}p^{rs}q(y)\right]dy=0
\end{array}
\right.
\label{eq:nearincompcell}
\end{equation}
for all $\phi\in H^1_\#(Y)^d$ and $q\in L^2(Y)$. We can then write
\[
u_1(\bz;x,y)=N^{rs}(\bz;x,y)\eps_{rs}(u_0(\bz;x)),\ \ \ p(\bz;x,y)=p^{rs}(\bz;x,y)\eps_{rs}(u_0(\bz;x)).
\]
Let $\mu^0$ be the fourth order tensor and $\lambda^0$ be the second order tensor defined by
\[
\mu^0_{ijrs}(\bz;x)=2\int_Y\mu(\bz;x,y)(e^{rs}_{ij}+\eps_{yij}(N^{rs}(\bz;x,y)))dy,\ \ 
\lambda^0_{rs}(\bz;x)=\int_Yp^{rs}(\bz;x,y)dy.
\]
The homogenized elastic tensor is 
$
a^0_{ijrs}(\bz;x)=\mu^0_{ijrs}(\bz;x)+\delta_{ij}\lambda^0_{rs}(\bz;x)
$
($a^0$ may not be isotropic). 
The homogenized equation, in the variational form, is
\[
\int_D[\mu^0_{ijrs}\eps_{rs}(u_0(\bz;x))\eps_{ij}(\phi(x))+\lambda^0_{rs}(\bz;x)\eps_{rs}(u^0_{rs}(\bz;x))\delta_{ij}\eps_{ij}(\phi(x))]dx=\int_Df(x)\cdot\phi(x)dx.
\]
We note that as $p^{rs}(\bz;x,y)=\lambda(\bz;x,y)(e_{ii}^{rs}+{\rm div}_yN^{rs}(\bz;x,y))$, this formula of $a^0_{ijrs}$ is consistent with \eqref{eq:a0}. Let 
\[
u^{1\ep}(\bz;x)=u^0(\bz;x)+\ep N^{rs}(\bz;x,{x\over\ep})\eps_{rs}(u^0)(\bz;x)
\ \ 
\mbox{and}\ \  
p^{1\ep}(\bz;x)=p^{rs}(\bz;x,{x\over\ep})\eps_{rs}(u^0(\bz;x)).
\]
For each function $\phi\in V$, we have
\beqas
&&\int_D[2\mu^\ep(\bz;x)\eps(u^{1\ep}(\bz;x)):\eps(\phi)+{\rm div}\phi(x)p^{1\ep}(\bz;x)]dx\\
&&=\int_D\left[2\mu^\ep(\bz;x)\left(\eps_{ij}(u^0(\bz;x))+\eps_{yij}(N^{rs}(\bz;x,{x\over\ep}))\eps_{rs}(u^0(\bz;x))+\ep\eps_{ij}(N^{rs}(\bz;x,{x\over\ep}))\eps_{rs}(u^0(\bz;x))\right.\right.\\
&&\qquad\qquad\left.+\frac12\ep\left(N_i^{rs}(\bz;x,\by){\partial\eps_{rs}(u_0(\bz;x))\over\partial x_j}+N_j^{rs}(\bz;x,\by){\partial\eps_{rs}(u_0(\bz;x))\over\partial x_i}\right)\right)\eps_{ij}(\phi)\\
&&\qquad\qquad+\left.\delta_{ij}p^{rs}(\bz;x,{x\over\ep})\eps_{rs}(u_0(\bz;x))\eps_{ij}(\phi(x))\right]dx\\
&&=\int_D(\mu^0_{ijrs}(\bz;x)+\delta_{ij}\lambda^0_{rs}(\bz;x))\eps_{rs}(u^0(\bz;x))\eps_{ij}(\phi(x))dx+
\int_Dg_{ijkl}(\bz;x,{x\over\ep})\eps_{kl}(u^0(x))\eps_{ij}(\phi(x))\\
&&\qquad\qquad+\ep I_{ij}\eps_{ij}(\phi)
\eeqas
where
\be
g_{ijkl}(\bz;x,y)=\mu(\bz;x,y)(\delta_{ik}\delta_{jl}+\delta_{il}\delta_{jk})+2\mu(\bz;x,y)\eps_{yij}(N^{kl}(\bz;x,y))+\delta_{ij}p^{kl}(\bz;x,y)-\mu^0_{ijkl}(\bz;x)-\delta_{ij}\lambda_{kl}^0(\bz;x)
\label{eq:gincomp}
\ee
and 
\[
I_{ij}=\eps_{ij}(N^{rs}(\bz;x,{x\over\ep}))\eps_{rs}(u^0(\bz;x))+\frac12\left(N_i^{rs}(\bz;x,\by){\partial\eps_{rs}(u_0(\bz;x))\over\partial x_j}+N_j^{rs}(\bz;x,\by){\partial\eps_{rs}(u_0(\bz;x))\over\partial x_i}\right)
\]
From \eqref{eq:nearincompcell}, we deduce that
\[
\int_Yg_{ijkl}(\bz;x,y)dy=0,\ \ \mbox{and}\ \ {\partial\over\partial y_j}g_{ijkl}(\bz;x,y)=0.
\]
From the result in \cite{JKO} page 7, there are functions 
\be
\alpha^{ikl}_{jr}(\bz;x,\cdot)\in H^1_\#(Y)\mbox{\ \  with\ \ } \alpha^{ikl}_{jr}=\alpha^{ikl}_{rj}\mbox{\ \  and\ \ } 
g_{ijkl}(\bz;x,y)={\partial\over\partial y_r}\alpha^{ikl}_{jr}(\bz;x,y).
\label{eq:alphatensor}
\ee
Since $g_{ijkl}=g_{jikl}$, 
\beqas
&&\int_Dg_{ijkl}(\bz;x,{x\over\ep})\eps_{kl}(u^0(\bz;x))\eps_{ij}(\phi(x))dx=\\
&&-\ep\int_D\alpha^{ikl}_{jr}(\bz;x,{x\over\ep})\eps_{kl}(u^0(\bz;x)){\partial^2\phi_i\over\partial x_j\partial x_r}(x)dx -\ep\int_D r_{ij}^\ep(\bz;x)\eps_{ij}(\phi(x))dx
=\ep\int_D r_{ij}^\ep(\bz;x)\eps_{ij}(\phi(x))dx
\eeqas
where
\[
r_{ij}^\ep(\bz;x)=-{\partial\alpha^{ikl}_{jr}\over\partial x_r}(\bz;x,{x\over\ep})\eps_{kl}(u^0)(\bz;x)-\alpha^{ikl}_{jr}(\bz;x,{x\over\ep}){\partial\over\partial x_r}\eps_{kl}(u^0(\bz;x)).
\]
Thus there is a constant $c=c(||\alpha||,||N||,||u^0||)$ so that
\beqas
&&\Bigg|\int_D[2\mu^\ep(\bz;x)\eps(u^{1\ep}(\bz;x)):\eps(\phi)+{\rm div}\phi(x)p^{1\ep}(\bz;x)]dx-\\
&&\qquad\qquad\int_D(\mu^0_{ijrs}(\bz;x)+\delta_{ij}\lambda^0_{rs}(\bz;x))\eps_{rs}(u^0(\bz;x))\eps_{ij}(\phi(x))dx\Bigg|\le c\ep\|\phi\|_{V}.
\eeqas
As 
\beqas
&&\int_D(\mu^0_{ijrs}(\bz;x)+\delta_{ij}\lambda^0_{rs}(\bz;x))\eps_{rs}(u^0(\bz;x))\eps_{ij}(\phi(x))dx=\\
&&\int_D[2\mu^\ep(\bz;x)\eps(u^{\ep}(\bz;x)):\eps(\phi)+{\rm div}\phi(x)p^{\ep}(\bz;x)]dx
=\int_Df(x)\cdot \phi(x)dx,
\eeqas
we have
\[
\Bigg|\int_D[2\mu^\ep(\bz;x)\eps(u^{1\ep}(\bz;x)-u^{\ep}(\bz;x))):\eps(\phi)+{\rm div}\phi(x)(p^{1\ep}(\bz;x)-p^\ep(\bz;x))]dx\Bigg|\le c\ep\|\phi\|_V.
\]
Let $\tau^\ep\in {\cal D}(D)$ be such that $\tau^\ep(x)=1$ outside an $\ep$ neighbourhood $\tilde D^\ep\subset D$ of $\partial D$ and $\sup_{x\in D}\ep|\nabla\tau^\ep(x)|<c$ for all $\ep$. Let
\[
w^{1,\ep}(\bz;x)=u^0(\bz;x)+\ep\tau^\ep(x)N^{rs}(\bz;x,{x\over\ep})\ep_{rs}(u^0)(\bz;x).
\]
We have
\beqas
{\partial\over\partial x_j}(u^{1\ep}-w^{1\ep})_i(\bz;x)=-\ep{\partial\tau^\ep\over\partial x_j}(x)N^{rs}_i(\bz;x,{x\over\ep})\eps_{rs}(u^0)(\bz;x)+\ep(1-\tau^\ep(x)){\partial N^{rs}_i\over\partial x_j}(\bz;x,{x\over\ep})\eps_{rs}(u^0)(\bz;x)+\\(1-\tau^\ep(x)){\partial N^{rs}_i\over\partial y_j}(\bz;x,{x\over\ep})\eps_{rs}(u^0)(\bz;x)
+\ep(1-\tau^\ep(x))N^{rs}_i(\bz;x,{x\over\ep}){\partial\over\partial x_j}\eps_{rs}(u^0)(\bz;x).
\eeqas
We have that  $\|\phi\|^2_{L^2(\tilde D^\ep)}\le c\ep^2\|\phi\|_{H^1(D)}^2+c\ep\|\phi\|_{L^2(\partial D)}^2\le c\ep\|\phi\|_{H^1(D)}^2$ for all $\phi\in C^\infty(D)$ and therefore for all $\phi\in H^1(D)$. 
As $u^0\in L^\infty(U;H^2(D)^d)$ so $\eps_{rs}(u^0)\in L^\infty(U;H^1(D))$. Together with $N^{rs}\in L^\infty(U;C^1(\bar D,C^1(\bar Y)))^d$,
there is a constant $c=c(||N||,||u^0||)$ such that  for all $i=1,\ldots,d$
\be
\|u^{1\ep}_i(\bz)-w^{1\ep}_i(\bz)\|_{H^1(D)}\le c(||N||,||u^0||)\ep^{1/2},\ \ \forall\,\bz\in U.
\label{eq:u1w1incomp}
\ee
Therefore
\[
\Bigg|2\int_D\mu^\ep(\bz;x)\eps(w^{1\ep}(\bz;x)-u^{1\ep}(\bz;x))):\eps(\phi)dx\Bigg|\le c(||N||,||u^0||)\ep^{1/2}\|\phi\|_V.
\]
Thus
\begin{eqnarray}
\Bigg|\int_D[2\mu^\ep(\bz;x)\eps(w^{1\ep}(\bz;x)-u^{\ep}(\bz;x))):\eps(\phi)+{\rm div}\phi(x)(p^{1\ep}(\bz;x)-p^\ep(\bz;x))]dx\Bigg|\le\nonumber\\
(c(||\alpha||,||N||,||u^0||)\ep+ c(||N||,||u^0||))\ep^{1/2}\|\phi\|_V.
\label{eq:incompw1u}
\end{eqnarray}
We note that
\beqas
{\rm div}u^{1\ep}(\bz;x)&=&{\rm div}u^0(\bz;x)+{\rm div}_yN^{rs}(\bz;x,{x\over\ep})\eps_{rs}(u^0(\bz;x))\\
&&\qquad\qquad\qquad+\ep N^{rs}(\bz;x,{x\over\ep})\cdot{\rm grad}\eps_{rs}(u^0(\bz;x))+\ep{\rm div}_xN^{rs}(\bz;x,{x\over\ep})\eps_{rs}(u^0(\bz;x))\\
&=&(e^{rs}_{ii}+{\rm div}_yN^{rs}(\bz;x,{x\over\ep}))\eps_{rs}(u^0(\bz;x))+\ep N^{rs}(\bz;x,{x\over\ep})\cdot{\rm grad}\eps_{rs}(u^0(\bz;x))\\
&&\qquad\qquad\qquad+\ep{\rm div}_xN^{rs}(\bz;x,{x\over\ep})\eps_{rs}(u^0(\bz;x))\\
&=&{1\over\lambda^\ep(\bz;x)}p^{1\ep}(x)+\ep N^{rs}(\bz;x,{x\over\ep})\cdot{\rm grad}\eps_{rs}(u^0(\bz;x))+\ep{\rm div}_xN^{rs}(\bz;x,{x\over\ep})\eps_{rs}(u^0(\bz;x)).
\eeqas
Therefore
\[
\Bigg|\int_D[{\rm div}u^{1\ep}(\bz;x)q(x)-{1\over\lambda^\ep(\bz;x)}p^{1\ep}(\bz;x)q(x)]dx\Bigg|\le c(||N||,||u^0||)\ep\|q\|_H,
\]
where the constant $c(||N||,||u^0||)$ does not depend on $\lambda_{\min}$ when it goes to $\infty$. 
From \eqref{eq:u1w1incomp}, we deduce 
\[
\Bigg|\int_D[{\rm div}w^{1\ep}(x)q(x)-{1\over\lambda^\ep(\bz;x)}p^{1\ep}(\bz;x)q(x)]dx\Bigg|\le c(||N||,||u^0||)\ep^{1/2}\|q\|_H,
\]
so 
\be
\Bigg|\int_D\left[{\rm div}(w^{1\ep}(\bz;x)-\ue(\bz;x))q(x)-{1\over\lambda^\ep(\bz;x)}(p^{1\ep}(\bz;x)-p^\ep(\bz;x))q(x)\right]dx\Bigg|\le c(\|N\|,\|u^0\|)\ep^{1/2}\|q\|_H.
\label{eq:incompw1p1}
\ee
From \eqref{eq:incompw1u} and \eqref{eq:incompw1p1}, 
\[
\|w^{1\ep}(\bz;x)-u^\ep(\bz;x)\|\le c(||\alpha||,||N||,||u^0||)\ep+c(||N||,||u^0||)\ep^{1/2}
\]
where the constants do not depend on $\lambda(\bz;x,y)$ when $\lambda_{\min}$ goes to $\infty$. From this we get the conclution. \eproof
\begin{remark}
The constant $c$ in the homogenization error \eqref{eq:homerror} depends also on $\|\alpha\|$, $\|N\|$ and $\|u_0\|$ (see the detailed proof in \cite{XHelasticity}). It also depends explicitly on $\sup_{x,\by}\lambda/\inf_{x,\by}\mu$. The constants in \eqref{eq:homerrornearlyincomp} does not depends on this ratio, so this error is better than \eqref{eq:homerror}. 
\end{remark}
For the required regularity of $u^0$ and $N^{rs}$, we make the following assumption
\begin{assumption}\label{assum:mumlambdam} The functions $\mu_m$ and $\lambda_m$ in \eqref{eq:mulambda} belong to $C^2(D,C^2_\#(Y))^{d^4}$ such that
\[
\sum_{m=1}^\infty\|\mu_m\|_{C^2(D,C^2(Y))^{d^4}}+\|\lambda_m\|_{C^2(D,C^2(Y))^{d^4}}<\infty.
\]
\end{assumption}
Under this assumption, we have that $N^{rs}\in L^\infty(U;C^1(\bar D,C^1(\bar Y)))^d$ and if $\partial D$ belongs to the $C^1$ class and the forcing $f\in L^2(D)^d$, then  $u^0\in L^\infty(U;H^2(D)^d)$. The proof uses elliptic regularity results in \cite{Mclean} and is presented in \cite{XHelasticity}. We then have the following result.
\begin{proposition}
Under Assumption \ref{assum:mumlambdam}, if the boundary $\partial D$ belongs to the $C^1$ class and $f\in L^2(D)^d$, then there are constants $c_1=c_1(||\alpha||,||N||,||u^0||)$ and $c_2=c_2(||N||,||u^0||)$ such that
\[
\left\|\nabla\ue(\bz;\cdot)-\left[\nabla u^0(\bz;\cdot)+\nabla_yu^1\left(\bz;\cdot,{\cdot\over\ep}\right)\right]\right\|_{L^2(U,\rho;L^2(D)^{d\times d})}\le c_1\ep+c_2\ep^{1/2}.
\]
\end{proposition}
From this, we have the following approximation result.
\begin{theorem}
Assume that the boundary $\partial D$ belongs to the $C^1$ class and $f\in L^2(D)^d$.  
Let $\theta>0$ be a constant. If Assumptions \ref{assum:psummabilitygammadelta} and \ref{assum:mumlambdam} hold with $\kappa<{\sqrt{3}\over 1+\theta}$ then there is $\vartheta>0$ depending on $\theta,\mu$ such that the solution $\bu$ of \eqref{eq:b3} and \eqref{eq:b4} and the solution $\bu_\Lambda$ of problem \eqref{eq:B3Lprob} and \eqref{eq:B4Lprob} satisfy
\[
\|\nabla\ue-[\nabla u^0_{\Lambda_N}+{\cal U}^\ep(\nabla_yu^1_{\Lambda_N})\|_{L^2(U,\rho;L^2(D)^{d\times d})}\le c_1\ep+c_2\ep^{1/2}+c_3N^{-s}
\]
where $s=1/p-1/2$. 
The constant $c_1$ depends on $\|\alpha\|$, $\|N\|$ and $\|\bu^0\|$, the constant $c_2$ depends on $\|N\|$ and $\|\bu^0\|$, the constant $c_3$ depends on ${\bar\lambda_{\max}\over\bar\lambda_{\min}}$, $\|\{\|\bu_\nu\|_\bV\}\|_{\ell^p(\cF)}$ and $\|\{\|p_\nu\|_\bH\}\|_{\ell^p(\cF)}$. 

\end{theorem}

\subsection{Multiscale problems}
We summarize briefly the derivation of the homogenized equation for the multiscale case. Details can be found in \cite{XHelasticity}. Let $a^n(x,\by_n)=a(x,\by)$. For $m=1,\ldots,n-1$, the $m$th level homogenized coefficient $a^m(\bz;x,\by_m)$ is defined recursively as follows. 
Let $N_{m+1}^{rs}(\bz;x,\by_{m+1})\in V_{m+1}$ be the solution of the cell problem 
\be
\int_D\int_{\bY_{m+1}}a^{m+1}(\bz;x,\by_{m+1})(e^{rs}+\epsilon_{y_{m+1}}(N_{m+1}^{rs})):\epsilon_{y_{m+1}}(\phi)d\by_{m+1}dx=0
\label{eq:cell}
\ee
for all $\phi\in V_{m+1}$. 
The $m$th level homogenized elastic moduli $a^m(\bz;x,\by_m)$  is
\be
a^{m}_{ijkl}(\bz;x,\by_{m})=\int_{Y_{m+1}}a^{m+1}(\bz;x,\by_{m+1})(e^{kl}+\epsilon_{y_{m+1}}(N_{m+1}^{kl})):(e^{ij}+\epsilon_{y_{m+1}}(N_{m+1}^{ij}))dy_{m+1};
\label{eq:am}
\ee
$a^0(\bz;x)$ is the homogenized coefficient. The homogenized equation is
\be
-{\partial\over\partial x_j}(a_{ijkl}^0\epsilon_{kl}(u_0))=f_i.
\label{eq:homeqn}
\ee
We deduce the convergence for the  multiscale solution $\ue$ of the parametric multiscale problem \eqref{eq:paradisplacement} in this section. For problems with more than two scales, a homogenization rate of convergence similar to that in \eqref{eq:homerror} is not available. However, we can deduce a corrector for the case where $\ep_i/\ep_{i+1}$ is an integer for $i=1,\ldots,n-1$. We first define the operator ${\cal T}^\ep_n:L^1(D)\to L^1(D\times \bY)$ as
\[
{\cal T}_n^\ep(\phi)(x,\by)
=
\phi\Bigl(\ep_1\Bigl[{x\over\ep_1}\Bigr]
+
\ep_2\Bigl[{y_1\over\ep_2/\ep_1}\Bigr]
+
\ldots
+
\ep_n\Bigl[{y_{n-1}\over\ep_n/\ep_{n-1}}\Bigr]+\ep_ny_n\Bigr)
\]
where $\phi\in L^1(D)$ is understood to be zero outside $D$. For each $\bz\in U$, when $\ep\to 0$, the solution $\ue(\bz)$ and its $n+1$-scale convergence limit  $\bu=(u^0,\ldots,u^n)$ which is the solution of problem \eqref{eq:msparadisplacement} satisfies
\[
{\cal T}_n^\ep\left({\frac{\partial\ue_i}{\partial x_j}}\right)\wc{\partial u^0_{i}\over\partial x_j}+{\partial u^1_{i}\over\partial y_{1j}}+\cdots+{\partial u^n_{i}\over\partial y_{nj}}
\]
in $L^2(D\times\bY)$ so
\be
{\cal T}_n^\ep(\eps(\ue))\wc\eps(u^0)+\eps_{y_1}(u^1)+\cdots+\eps_{y_n}(u^n)\ \ \mbox{in}\ \ L^2(D\times\bY)^{d\times d}.
\label{eq:mshomprog}
\ee
Letting $D^{\ep}$ be the $2\ep$ neighbourhood of $D$. We have 
\be
\int_D\phi dx=\int_{D^{\ep}}\int_{Y_1}\cdots\int_{Y_n}{\cal T}_n^\ep(\phi)dy_n\cdots dy_1dx\ \ \ \forall\,\phi\in L^1(D).
\label{eq:Tepsuewc}
\ee

The proofs of \eqref{eq:mshomprog} and \eqref{eq:Tepsuewc} can be found in \cite{CDG}. To deduce an approximation of $\ue$ in $H^1(D)^d$  we define the operator ${\cal U}_n^\ep:L^2(D\times\bY)\to L^2(D)$ as 
\beqas
{\cal U}_n^\ep(\Phi)(x)
=
\int_{Y_1}\cdots\int_{Y_n}\Phi\Bigl(\ep_1\Bigl[{x\over\ep_1}\Bigr]
+
\ep_1 t_1,{\ep_2\over\ep_1}\Bigl[{\ep_1\over\ep_2}\Bigl\{{x\over\ep_1}\Bigr\}\Bigr]
+
{\ep_2\over\ep_1}t_2,\cdots,
\\
{\ep_n\over\ep_{n-1}}\Bigl[{\ep_{n-1}\over\ep_n}\Bigl\{{x\over\ep_{n-1}}\Bigr\}\Bigr]
+
{\ep_n\over\ep_{n-1}}t_n,\Bigl\{{x\over\ep_n}\Bigr\}\Bigr)dt_n\cdots dt_1
\eeqas
for all functions $\Phi\in L^1(D\times Y_1\times\cdots\times Y_n)$. We assume further regularity for the elastic moduli.
\begin{assumption} \label{assum:psimms} The fourth order tensors $\psi_m$ in \eqref{eq:randommoduli} belong to $C^2(D,C^2_\#(Y_1,\ldots,C^2_\#(Y_n)\dots))^{d^4}$, which we denote as $C^2(D,C^2_\#(\bY))^{d^4}$, such that
$
\sum_{m=1}^\infty\|\psi_m\|_{C^2(D,C^2_\#(\bY))^{d^4}}<\infty
$
and $\bar a\in C^2(D,C^2_\#(\bY))^{d^4} $. 
\end{assumption}
We then have the following regularity result for the solution $N^{rs}$ of the  cell problem \eqref{eq:cell}. 
\begin{lemma}\label{lem:Nrsregularityms}
Under Assumption \ref{assum:psimms}, $N^{rs}_i\in L^\infty(U;C^1(\bar D,C^1_\#(\bY_i)))$ $\forall\,r,s=1,\ldots,d$ and $i=1,\ldots,n$. 
\end{lemma}
The proof is a routine generalization of the proof of Lemma \ref{lem:Nrsregularity}. To deduce the correctors for $\ue(\bz)$, we employ the following result which is established in Xia and Hoang \cite{XHelasticityrandom} for non-parametric problems. 
\begin{lemma}\label{lem:l3}
Assume that $u^0\in L^\infty(U;H^2(D)^d)$ and $N^{rs}_i\in L^\infty(U;C^1(\bar D,C^1_\#(\bY_i))^d)$. Then
\be
\lim_{\ep\to 0}\sup_{\bz\in U}\int_D\left|\eps_{y_i}(u^i(\bz;x,{x\over\ep_1},\ldots,{x\over\ep_i}))-{\cal U}^\ep_n(\eps_{y_i}(u^i))(\bz;x)\right|^2dx=0
\label{eq:e}
\ee
and 
\be
\lim_{\ep\to 0}\sup_{\bz\in U}\int_D\left|\nabla_{y_i}u^i(\bz;x,{x\over\ep_1},\ldots,{x\over\ep_i})-{\cal U}^\ep_n(\nabla_{y_i}u^i)(\bz;x)\right|^2dx=0.
\label{eq:e1}
\ee
\end{lemma}
The proof is similar to that for the non-parametric case presented in \cite{XHelasticity}; we refer to \cite{XHelasticity} for details. 
We then have the following convergence result.
\begin{theorem}
Assume that $D$ is a $C^1$ domain and $f\in L^2(D)^d$. Under Assumption \ref{assum:psimms}, 
\[
\lim_{\ep\to 0\atop N\to\infty}\|\nabla\ue-{\cal U}_n^\ep(\nabla u^0_\Lambda+\nabla_{y_1}u^1_{\Lambda}+\ldots+\nabla_{y_n}u^n_\Lambda)\|_{L^2(U,\rho;L^2(D)^{d\times d})}=0
\]
where $\bu_\Lambda=(u^0_\Lambda,u^1_\Lambda,\ldots,u^n_\Lambda)$ is the solution of problem \eqref{eq:GBprob} corresponding to the best $N$ term set $\Lambda$.  
\end{theorem}
\bproof Following the proof of Lemma 6.6 of \cite{XHelasticity}, we have
\[
\lim_{\ep\to 0}\|\nabla\ue-{\cal U}^\ep_n(\nabla u^0+\ldots+\nabla_{y_n}u^n)\|_{L^2(U,\rho;L^2(D)^{d\times d})}=0.
\]
We get the conclusion by using
\[
\|{\cal U}^\ep_n(\nabla_{y_i}u^i-\nabla_{y_i}u^i_\Lambda)\|_{L^2(U,\rho;L^2(D)^{d\times d})}\le \|\nabla_{y_i}u^i-\nabla_{y_i}u^i_\Lambda\|_{L^2(U,\rho;L^2(D\times\bY)^{d\times d})}\to0\mbox{\ when\ }N\to\infty.
\]
\eproof

For the mixed problems, we have the following result.
\begin{theorem} Assume that the boundary $\partial D$ belongs to $C^1$ and $f\in L^2(D)^d$. Under Assumption \ref{assum:psimms}
\[
\lim_{\ep\to 0\atop N\to \infty}\|\nabla\ue-[\nabla u^0_{\Lambda_N}+{\cal U}^\ep_n(\nabla_{y_1}u^1_{\Lambda}+\ldots+\nabla_{y_n}u^n_\Lambda)]\|_{L^2(U,\rho;L^2(D)^{d\times d})}+\|\sigma^\ep-{\cal U}^\ep(\sigma_{\Lambda})\|_{L^2(U,\rho;L^2(D)^{d\times d})}=0.
\]
\end{theorem}
\bproof
From \eqref{eq:e}, we have 
$
\lim_{\ep\to 0}\|\eps(\ue)-{\cal U}^\ep_n(\bep(\bu))\|_{L^2(U,\rho;L^2(D)^{d\times d})}=0,
$
so\\
$
\lim_{\ep\to 0}\|\sigma^\ep-a^\ep{\cal U}^\ep_n(\bep(\bu))\|_{L^2(U,\rho;L^2(D)^{d\times d})}=0.
$
Since $a\in L^\infty(U;C^2(D,C^2_\#(\bY)))^{d^4}$, 
\[
\lim_{\ep\to 0}\sup_{\bz\in U}\sup_{x\in D}|a(\bz;x,{x\over\ep_1},\ldots,{x\over\ep_n})-{\cal U}^\ep_n(a)(x)|=0.
\]
Therefore 
$
\lim_{\ep\to 0}\|\sigma^\ep-{\cal U}^\ep_n(a\bep(\bu))\|_{L^2(U,\rho;L^2(D)^{d\times d})}=0,
$
i.e.
$
\lim_{\ep\to 0}\|\sigma^\ep-{\cal U}^\ep_n(\sigma)\|_{L^2(U,\rho;L^2(D)^{d\times d})}=0.
$
From Theorem \ref{thm:bestNterm12},
$
\lim_{N\to\infty}\|{\cal U}^\ep_n(\sigma)-{\cal U}^\ep_n(\sigma_{\Lambda})\|_{L^2(U,\rho;L^2(D)^{d\times d})}=0.
$
Therefore
\[
\lim_{\ep\to 0\atop N\to\infty}\|\nabla\ue-[\nabla u^0_{\Lambda_N}+{\cal U}^\ep_n(\nabla_{y_1}u^1_\Lambda+\ldots+\nabla_{y_n}u^n_\Lambda)]\|_{L^2(U,\rho;L^2(D)^{d\times d})}+\|\sigma^\ep-{\cal U}^\ep_n(\sigma_{\Lambda})\|_{L^2(U,\rho;L^2(D)^{d\times d})}=0.
\]
\eproof

\begin{remark} The purpose of considering mixed problems with a penalty term for nearly incompressible materials is to derive an approximation whose error is independent of the ratio of the Lam\'e constants. For the general multiscale problems, we do not have an explicit error for the corrector so we do not consider the nearly incompressible problems separately in this section.
\end{remark}
\bibliographystyle{plain}
\bibliography{randommulti}

\end{document}